%% file: main.tex
\documentclass[10pt,conference]{IEEEtran}
\usepackage{amsmath}
\usepackage{mathtools}
\usepackage{epsfig}
\usepackage{amssymb}
\usepackage{amsfonts}
\usepackage{hyperref}
\usepackage{url}
\usepackage[hyphenbreaks]{breakurl}

\author{
	\IEEEauthorblockN{
		Manda Winlaw\IEEEauthorrefmark{1},
		Michael B Hynes\IEEEauthorrefmark{2},
		Anthony Caterini\IEEEauthorrefmark{3},
		Hans De Sterck\IEEEauthorrefmark{4}$^{,1}$
	}
	\IEEEauthorblockA{
		Department of Applied Mathematics\\
		University of Waterloo, Canada\\
		Email: 
		\IEEEauthorrefmark{1}mwinlaw@uwaterloo.ca,
		\IEEEauthorrefmark{2}mbhynes@uwaterloo.ca,
		\IEEEauthorrefmark{3}alcaterini@uwaterloo.ca,
		\IEEEauthorrefmark{4}hdesterck@uwaterloo.ca
	}
}

\newcommand{\insertKeywords}{
	\begin{IEEEkeywords}
		Recommendation systems,
		collaborative filtering,
		parallel optimization algorithms,
		matrix factorization,
		Apache Spark,
		Big Data,
		scalable methods.
	\end{IEEEkeywords}
}
\newcommand{\insertBibliography}{
	{
		\balance
		\bibliographystyle{IEEEtran}
		\bibliography{AcceleratedALS} 
	}
}
\title{
	Algorithmic Acceleration of Parallel ALS for Collaborative Filtering: Speeding up Distributed Big Data Recommendation in Spark
	\vspace{-0.3cm}
}

\usepackage[ruled,vlined]{algorithm2e}
\usepackage{epstopdf}
\usepackage{balance}
\usepackage{color}

\graphicspath{{fig/}}

\newcommand{\eqn}[1]{(\ref{Equation:#1})}
\newcommand{\fig}[1]{Fig.~\ref{Figure:#1}}
\newcommand{\floor}[1]{\left\lfloor{#1} \right\rfloor }
\newcommand{\ceil}[1]{\left\lceil{#1} \right\rceil }
\newcommand{\norm}[1]{\|#1\|}
\newcommand{\uvec}{\mathbf{u}_i}
\newcommand{\mvec}{\mathbf{m}_j}
\newcommand{\vvec}{\mathbf{v}}
\newcommand{\gvec}{\mathbf{g}}
\newcommand{\gpvec}{\mathbf{\bar{g}}}
\newcommand{\pvec}{\mathbf{p}}

\newcommand{\xvec}{\mathbf{x}}
\newcommand{\yvec}{\mathbf{y}}
\newcommand{\zvec}{\mathbf{z}}

\newcommand{\xpvec}{\mathbf{\bar{x}}}
\newcommand{\Imat}{\mathbf{I}}
\newcommand{\Umat}{\mathbf{U}}
\newcommand{\Mmat}{\mathbf{M}}
\newcommand{\Rmat}{\mathbf{R}}
\newcommand{\Amat}{\mathbf{A}}

\newcommand{\Iset}{\mathcal{I}}
\newcommand{\betap}{\bar{\beta}_{k+1}}
\newcommand{\lineage}{\mathcal{L}}
\newcommand{\rdd}[1]{\mathtt{\small{RDD}}\left[  #1 \right] }
\newcommand{\Set}[1]{\left\{ #1 \right\}}

\newcommand{\key}{\mathcal{K}}
\newcommand{\A}{\mathcal{A}}
\newcommand{\B}{\mathcal{B}}
\newcommand{\order}[1]{\mathcal{O}( #1 )}

\begin{document}
\maketitle

\begin{abstract}
Collaborative filtering algorithms are important building blocks in many practical recommendation systems.
For example, many large-scale data processing environments include collaborative filtering models for which the Alternating Least Squares (ALS) algorithm is used to compute latent factor matrix decompositions.
In this paper, we propose an approach to accelerate the convergence of parallel ALS-based optimization methods for collaborative filtering using a nonlinear conjugate gradient (NCG) wrapper around the ALS iterations. 
We also provide a parallel implementation of the accelerated
ALS-NCG algorithm in the Apache Spark distributed data
processing environment, and an efficient line search
technique as part of the ALS-NCG implementation that requires only one pass over the data on distributed datasets.
In serial numerical experiments on a linux workstation and
parallel numerical experiments on a 16 node cluster with 256
computing cores, we demonstrate that the combined ALS-NCG
method requires many fewer iterations and less time than
standalone ALS to reach movie rankings with high accuracy on 
the MovieLens 20M dataset. 
In parallel, ALS-NCG can achieve an acceleration factor of 4
or greater in clock time when an accurate solution is
desired; furthermore, the acceleration factor increases
as greater numerical precision is required in the
solution.
In addition, the NCG acceleration mechanism is efficient
in parallel and scales linearly with problem size on
synthetic datasets with up to nearly 1 billion ratings.
The acceleration mechanism is general and may also be
applicable to other optimization methods for collaborative
filtering.
\end{abstract}

\insertKeywords
\section{Introduction and Background}

\footnotetext[1]{Currently at Monash University, School of Mathematical Sciences, Melbourne, Australia}

Recommendation systems are designed to analyze available user data to
recommend items such as movies, music, or other goods to consumers, and
have become an increasingly important part of most successful online businesses.
One strategy for building recommendation systems is known as collaborative filtering,
whereby items are recommended to users by collecting preferences or taste
information from many users
(see, e.g., \cite{Sarwar:2001,koren2008factorization}).
Collaborative filtering methods provide the basis for
many recommendation systems \cite{Bobadilla:2013} and have been used by online
businesses such as Amazon \cite{Linden:2003}, Netflix \cite{Bell:2007a}, and Spotify
\cite{Johnson:2014}.

An important class of collaborative filtering methods are latent factor models, for which low-rank matrix factorizations are
often used and have
repeatedly demonstrated better accuracy than other
methods such as nearest neighbor models and restricted Boltzmann machines
\cite{Bell:2007a,Dror:2012}.
A low-rank latent factor model associates with each user and with each item a vector of rank $n_f$, for which each component measures the tendency of the user or item towards a certain \emph{factor} or \emph{feature}. 
In the context of movies, a latent feature may represent the
style or genre (i.e.~drama, comedy, or romance), and the magnitude of a given component of a
user's feature vector is proportional to the user's proclivity for that feature
(or, for an item's feature vector, to the degree to which that feature is manifested by the item).
A low-rank matrix factorization procedure takes as its input
the user-item ratings matrix $\mathbf{R}$, in which each
entry is a numerical rating of an item by a user and for which
typically very few entries are known.
The procedure then
determines the low-rank user
($\mathbf{U}$) and item ($\mathbf{M}$) matrices,
where a column in these matrices represents a latent feature vector for a single user or item, respectively, and $\Rmat \approx \mathbf{U}^T\mathbf{M}$ for the known values in $\mathbf{R}$.
Once the user and item matrices are computed, they are used to build the
recommendation system and predict the unknown ratings.
Computing user and
item matrices is the first step in building a variety of recommendation systems,
so it is important to compute the factorization of $\Rmat$
quickly.

The matrix factorization problem is closely related to the
singular value decomposition (SVD), but the SVD of a matrix
with missing values is undefined. However, since $\mathbf{R}
\approx \mathbf{U}^T\mathbf{M}$, one way to find the user
and item matrices is by minimizing the squared difference
between the approximated and actual value of the known ratings in $\Rmat$.  Minimizing this difference is typically done by one of two algorithms: stochastic gradient descent (SGD) or alternating least squares (ALS) \cite{Koren:2009,Funk:2006}.  
ALS can be easily parallelized and can efficiently handle models that incorporate implicit data (e.g. the frequency of a user's mouse-clicks or time spent on a website) \cite{Hu:2008}, but it is well-known that ALS can require a large number of iterations to converge. Thus, in this paper, we propose an approach to significantly accelerate the convergence of the ALS algorithm for computing the user and item matrices.  We use a nonlinear optimization algorithm, specifically the nonlinear conjugate gradient (NCG) algorithm \cite{Nocedal:2006}, as a wrapper around ALS to significantly accelerate the convergence of ALS, and thus refer to this combined algorithm as ALS-NCG.  Alternatively, the algorithm can be viewed as a nonlinearly preconditioned NCG method with ALS as the nonlinear preconditioner \cite{DeSterck:2015}.
Our approach for accelerating the ALS algorithm using the NCG algorithm can be situated in the context of recent research activity on nonlinear preconditioning for nonlinear iterative solvers \cite{PETSC:2015,DeSterck:2012,DeSterck:2015,Fang:2009,Walker:2011}.  Some of the ideas date back as far as the 1960s \cite{Anderson:1965,Concus:1977}, but they are not well-known and remain under-explored experimentally and theoretically \cite{PETSC:2015}.  

Parallel versions of collaborative filtering and recommendation are of great interest in the era of big data
\cite{gemulla2011large,teflioudi2012distributed,yu2012scalable}.
For example, the Apache Spark data processing environment \cite{zaharia2012resilient}
contains a parallel implementation of ALS for the
collaborative filtering model of \cite{Zhou:2008,Koren:2009}.
In \cite{gemulla2011large} an advanced distributed SGD method is described
in Hadoop environments, followed by work in \cite{teflioudi2012distributed} that considers algorithms based
on ALS and SGD in environments that use the Message Passing Interface (MPI).
Scalable coordinate descent approaches for collaborative filtering are proposed in \cite{yu2012scalable},
also using the MPI framework.
In addition to producing algorithmic advances, these papers
have shown that the relative performance of the methods considered is often
strongly influenced by the performance characteristics of
the parallel computing paradigm that is used to implement
the algorithms (e.g., Hadoop or MPI).
In this paper we implement our proposed ALS-NCG algorithm in
parallel using the Apache Spark framework, which is a large-scale distributed data processing environment that builds on the principles of scalability and fault tolerance that are instrumental in the success of Hadoop and MapReduce. 
However, Spark adds crucial new capabilities in terms of in-memory computing and data persistence for iterative methods, and makes it possible to implement more elaborate algorithms with significantly better performance than is
feasible in Hadoop. 
As such, Spark is already being used extensively for advanced big data analytics in the commercial setting \cite{Johnson2014}.
The parallel Spark implementation of ALS for the collaborative filtering model of \cite{Zhou:2008,Koren:2009}
 forms the starting point for applying the parallel acceleration methods proposed in this paper.

Our contributions in this paper are as follows. The specific optimization problem is formulated in Section \ref{Section:PD}, and the accelerated ALS-NCG algorithm  is developed in Section \ref{Section:OA}. In Section \ref{Section:NR}, we study the convergence enhancements of ALS-NCG compared to standalone ALS in small serial tests using subsets of the MovieLens 20M dataset \cite{MovieLens:2015}.
Section \ref{Section:Parallel} describes our parallel
implementation in Spark\footnote{Source code available:
	\url{https://github.com/mbhynes/als-ncg}}, 
and Section \ref{Section:sparkExperiments} contains results from parallel performance tests on a high-end computing cluster with 16 nodes and 256 cores, using both the full MovieLens 20M dataset and a large synthetic dataset sampled from the MovieLens 20M data with up to 6 million users and 800 million ratings. 
We find that ALS-NCG converges significantly faster than ALS in both the serial and distributed Spark settings, demonstrating the overall speedup provided by our algorithmic acceleration.

\section{Problem Description}\label{Section:PD}
The acceleration approach we propose in this paper is applicable to a broad class of optimization methods and collaborative filtering models. For definiteness, we choose a specific latent factor model, the matrix factorization model from \cite{Koren:2009} and \cite{Zhou:2008}, and a specific optimization method, ALS \cite{Zhou:2008}. Given the data we use, the model is presented in terms of users and movies instead of the more generic users and items framework.

Let the matrix of user-movie rankings be represented by $\mathbf{R} = \{r_{ij}\}_{n_u\times n_m}$ where $r_{ij}$ is the rating given to movie $j$ by user $i$, $n_u$ is the number of users, and $n_m$ is the number of items.
Note that for any user $i$ and movie $j$, the value of $r_{ij}$ is either a real number or is missing, and in practice very few values are known.
For example, the MovieLens 20M dataset \cite{MovieLens:2015} with 138,493 users and 27,278 movies contains only 20 million rankings, accounting for less than 1$\%$ of the total possible rankings.  In the low-rank factorization of $\mathbf{R}$ with $n_f$ factors, each user $i$ is associated with a vector $\mathbf{u}_{i} \in \mathbb{R}^{n_f}$ ($i = 1,\ldots,n_u$), and each movie $j$ is associated with a vector $\mathbf{m}_{j} \in \mathbb{R}^{n_f}$ ($j = 1,\ldots,n_m$). The elements of $\mathbf{m}_{j}$ measure the degree that movie $j$ possesses each factor or feature, and the elements of $\mathbf{u}_{i}$ similarly measures the affinity of user $i$ for each factor or feature.
The dot product $\mathbf{u}_i^T\mathbf{m}_j$ thus captures the interaction between user $i$ and movie $j$, approximating user $i$'s rating of movie $j$ as $r_{ij} \approx \mathbf{u}_i^T\mathbf{m}_j$.
Denoting $\mathbf{U} = \left[\mathbf{u}_i\right] \in
\mathbb{R}^{n_f \times n_u}$ as the user feature matrix and
$\mathbf{M} = \left[\mathbf{m}_j\right] \in \mathbb{R}^{n_f \times n_m}$ as the movie feature matrix, our goal is to determine $\Umat$ and $\Mmat$ such that $\mathbf{R} \approx \mathbf{U}^T\mathbf{M}$ by minimizing the following squared loss function:
\begin{equation}
\label{Equation:ActualLoss}
\begin{split}
\mathcal{L}_{\lambda}(\mathbf{R},\mathbf{U},\mathbf{M})  &= \sum_{(i,j)\in \mathcal{I}} (r_{ij} - \mathbf{u}_i^T\mathbf{m}_j)^2 + \\
&\lambda\Big(\sum_i n_{u_i}\|\mathbf{u}_i\|^2 + \sum_j n_{m_j}\|\mathbf{m}_j\|^2\Big),
\end{split}
\end{equation} 
where $\mathcal{I}$ is the index set of known $r_{ij}$ in $\Rmat$, $n_{u_i}$ denotes the number of ratings by user $i$, and $n_{m_j}$ is the number of ratings of movie $j$.  The term $\lambda(\sum_i n_{u_i}\|\mathbf{u}_i\|^2 + \sum_j n_{m_j}\|\mathbf{m}_j\|^2)$ is a Tikhonov regularization \cite{Tikhonov:77} term commonly included in the loss function to prevent overfitting.  The full optimization problem can be stated as 
\begin{equation}
\label{Equation:OptProblem}
\min_{\mathbf{U},\mathbf{M}} \mathcal{L}_{\lambda}(\mathbf{R},\mathbf{U},\mathbf{M}).
\end{equation}

\section{Accelerating ALS Convergence by NCG}\label{Section:OA}

\subsection{Alternating Least Squares Algorithm}
The optimization problem in (\ref{Equation:OptProblem}) is not convex.
However, if we fix one of the unknowns, either $\mathbf{U}$ or $\mathbf{M}$, then the optimization problem becomes quadratic and we can solve for the remaining unknown as a least squares problem. Doing this in an alternating fashion is the central idea behind ALS. 

Consider the first step of the ALS algorithm in which $\mathbf{M}$ is fixed. 
We can determine the least squares solution to \eqn{ActualLoss} for each $\uvec$ by setting all the components of the gradient of \eqn{ActualLoss} related to $\uvec$ to zero: $\frac{\partial\mathcal{L}_{\lambda}}{\partial u_{ki}} = 0 \;\; \forall \; i,k$ where $u_{ki}$ is an element of $\Umat$.
Expanding the terms in the derivative of \eqn{ActualLoss}, we have
\begin{equation*}
\sum_{j \in \mathcal{I}_i}2(\mathbf{u}_i^T\mathbf{m}_j - r_{ij})m_{kj} + 2\lambda n_{u_i}u_{ki} = 0 \;\; \forall \; i,k
\end{equation*}
\begin{equation*}
\Rightarrow  \;\; \sum_{j \in \mathcal{I}_i}m_{kj}\mathbf{u}_i^T\mathbf{m}_j + \lambda n_{u_i}u_{ki} = \sum_{j \in \mathcal{I}_i}m_{kj}r_{ij} \;\; \forall \; i,k,
\end{equation*}
where $m_{kj}$ is an element of $\mathbf{M}$. 
In vector form, the resultant linear system for any $\uvec$ is
\begin{equation*}
(\mathbf{M}_{\mathcal{I}_i}\mathbf{M}_{\mathcal{I}_i}^T  + \lambda n_{u_i}\mathbf{I})\mathbf{u}_i =  \mathbf{M}_{\mathcal{I}_i}\mathbf{R}^T(i,\mathcal{I}_i) \;\; \forall \; i,
\end{equation*}
where $\mathbf{I}$ is the $n_f \times n_f$ identity matrix, $\mathcal{I}_i$ is the index set of movies user $i$ has rated, and $\mathbf{M}_{\mathcal{I}_i}$ represents the sub-matrix of $\mathbf{M}$ where columns $j \in \mathcal{I}_i$ are selected.
Similarly, $\mathbf{R}(i,\mathcal{I}_i)$ is a row vector that represents the $i$th row of $\mathbf{R}$ with only the columns in $\mathcal{I}_i$ included. 
The explicit solution for $\uvec$ is then given by
\begin{equation}
\label{Equation:SolveU}
	\mathbf{u}_i = \mathbf{A}_i^{-1}\mathbf{v}_i \;\; \forall \; i,
\end{equation}
where $\mathbf{A}_i = \mathbf{M}_{\mathcal{I}_i}\mathbf{M}_{\mathcal{I}_i}^T  + \lambda n_{u_i}\mathbf{I}$, and $\mathbf{v}_i = \mathbf{M}_{\mathcal{I}_i}\mathbf{R}^T(i,\mathcal{I}_i)$. 
The analogous solution for the columns of $\mathbf{M}$ is found by fixing $\mathbf{U}$, where each $\mvec$ is given by 
\begin{equation}
	\mvec = \mathbf{A}_j^{-1}\mathbf{v}_j \; \forall \; j,
\label{Equation:SolveM}
\end{equation}
where $\mathbf{A}_j = \mathbf{U}_{\mathcal{I}_j}\mathbf{U}_{\mathcal{I}_j}^T  + \lambda n_{m_j}\mathbf{I}$, $\mathbf{v}_j = \mathbf{U}_{\mathcal{I}_j}\mathbf{R}(\mathcal{I}_j,j)$.  
Here, $\mathcal{I}_j$ is the index set of users that have rated movie $j$, $\mathbf{U}_{\mathcal{I}_j}$ represents the sub-matrix of $\mathbf{U}  \in \mathbb{R}^{n_f \times n_u}$ where columns $i \in \mathcal{I}_j$ are selected, and $\mathbf{R}(\mathcal{I}_j,j)$ is a column vector that represents the $j$th column of $\mathbf{R}$ with only the rows in $\mathcal{I}_j$ included.

Algorithm \ref{Algorithm:ALS} summarizes the ALS algorithm used to solve the optimization problem given in (\ref{Equation:OptProblem}). From \eqn{SolveU} and \eqn{SolveM}, we note that each of the columns of $\mathbf{U}$ and $\mathbf{M}$ may be computed independently; thus ALS may be easily implemented in parallel, as in \cite{Zhou:2008}.
However, the ALS algorithm can require many iterations for convergence, and we now propose an acceleration method for ALS that can be applied to both the parallel and serial versions.
\begin{algorithm}[b]
 \SetAlgoLined
 \SetKwInOut{Output}{Output}
 \Output{$\Umat,\Mmat$}
   Initialize $\mathbf{M}$ with random values;\\
 \While{Stopping criteria have not been satisfied}{
  \For{$i = 1,\ldots,n_u$}{
	$\mathbf{u}_i \leftarrow \mathbf{A}_i^{-1}\mathbf{v}_i$;
	}
	\For{$j = 1,\ldots,n_m$}{
	$\mvec \leftarrow \mathbf{A}_j^{-1}\mathbf{v}_j$;
	}
 }
 \caption{Alternating Least Squares (ALS)}
 \label{Algorithm:ALS}
\end{algorithm}
\subsection{Accelerated ALS Algorithm}
In this section, we develop the accelerated ALS-NCG algorithm to solve the collaborative filtering optimization problem minimizing (\ref{Equation:ActualLoss}).
In practice we found that the nonlinear conjugate gradient algorithm by itself was very slow to converge to a solution of (\ref{Equation:OptProblem}), and thus do not consider it as an alternative to the ALS algorithm.  Instead we propose using NCG to accelerate ALS, or, said differently, the combined algorithm uses ALS as a nonlinear preconditioner for NCG \cite{DeSterck:2015}.

The standard NCG algorithm is a line search algorithm in continuous optimization that is an extension of the CG algorithm for linear systems. 
While the linear CG algorithm is specifically designed to minimize the convex quadratic function $\phi(\mathbf{x}) = \frac{1}{2}\mathbf{x}^T\mathbf{A}\mathbf{x} - \mathbf{b}^T\mathbf{x}$, where $\mathbf{A} \in \mathbb{R}^{n\times n}$ is a symmetric positive definite matrix, the NCG algorithm can be applied to general constrained optimization problems of the form $\min_{\textbf{x} \in \mathbb{R}^n}f(\textbf{x}).$
Here, the minimization problem is (\ref{Equation:OptProblem}), where the matrices $\mathbf{U}$ and $\mathbf{M}$ are found by defining the vector $\xvec \in \mathbb{R}^{n_f\times(n_u+n_m)}$ as
\begin{equation}
	\label{Equation:x}
	\mathbf{x}^T = \left[ \mathbf{u}^T_1 \; \mathbf{u}^T_2 \ldots \mathbf{u}^T_{n_u} \; \mathbf{m}^T_1 \; \mathbf{m}^T_2 \ldots \mathbf{m}^T_{n_m}\right] 
\end{equation}
with function $f(\xvec) = \mathcal{L}_{\lambda}$ as in \eqn{ActualLoss}.

The NCG algorithm generates a sequence of iterates $\mathbf{x}_i$, $i \geq 1$, from the initial guess $\mathbf{x}_0$ using the recurrence relation
\begin{equation*}
\label{Equation:NCG:StepUpdate}
\mathbf{x}_{k+1} = \mathbf{x}_k + \alpha_k \mathbf{p}_k.
\end{equation*}
The parameter $\alpha_k > 0$ is the step length determined by a line search along search direction $\mathbf{p}_k$, which is generated by the following rule:
\begin{equation}
\label{Equation:NCG:SearchDirectionUpdate}
\mathbf{p}_{k+1} = -\mathbf{g}_{k+1} + \beta_{k+1} \mathbf{p}_k, \;\; \mathbf{p}_0 = -\mathbf{g}_0,
\end{equation}
where $\beta_{k+1}$ is the update parameter and $\mathbf{g}_k = \nabla f(\mathbf{x}_k)$ is the gradient of $f(\xvec)$ evaluated at $\mathbf{x}_k$.
The update parameter $\beta_{k+1}$ can take on various different forms.
In this paper, we use the variant of $\beta_{k+1}$ developed by Polak and Ribi\`{e}re \cite{Polak:1969}:
\begin{equation}
\label{Equation:NCG:PRBeta}
\beta_{k+1} =
\frac{\mathbf{g}^T_{k+1}(\mathbf{g}_{k+1}-\mathbf{g}_k)}{\mathbf{g}_{k}^T\mathbf{g}_k}.
\end{equation}
Note that if a convex quadratic function is optimized using
the NCG algorithm with an exact line search, then \eqn{NCG:PRBeta} reduces to the same $\beta_{k+1}$ as in the original CG algorithm for linear systems \cite{Nocedal:2006}.

The preconditioning of the NCG algorithm by the ALS algorithm modifies the expressions for $\beta_{k+1}$ and $\pvec_{k+1}$ as follows, and is summarized in Algorithm \ref{Algorithm:ALS-NCG}.
Let $\overline{\mathbf{x}}_k$ be the iterate generated from one iteration of the ALS algorithm applied to $\xvec_k$, $\overline{\mathbf{x}}_k = P(\mathbf{x}_k)$,
where $P$ represents one iteration of ALS.
This iterate is incorporated into the NCG algorithm by defining the preconditioned gradient direction generated by ALS as
\begin{equation*}
\overline{\mathbf{g}}_k = {\mathbf{x}}_k - \overline{\mathbf{x}}_k = {\mathbf{x}}_k - P({\mathbf{x}}_k),
\end{equation*}
and replacing ${\mathbf{g}}_k$ with $\overline{\mathbf{g}}_k$ in \eqn{NCG:SearchDirectionUpdate}.
Note that $-\overline{\mathbf{g}}_k$ is expected to be a descent direction that is an improvement compared to the steepest descent direction $-\mathbf{g}_k$.
The update parameter $\beta_{k+1}$ is redefined as $\overline{\beta}_{k+1}$ with the form 
\begin{equation}
\label{Equation:AALS:PRBeta}
\overline{\beta}_{k+1} =
\frac{\overline{\gvec}_{k+1}^T(\gvec_{k+1} - \gvec_k
)}{{\overline{\mathbf{g}}}_{k}^T \mathbf{g}_{k}},
\end{equation}
where \eqn{AALS:PRBeta} is similar to \eqn{NCG:PRBeta}, however not every instance of ${\mathbf{g}}_{k}$ has been replaced with $\overline{\mathbf{g}}_{k}$.
The specific form for $\overline{\beta}_{k+1}$
is chosen because if Algorithm \ref{Algorithm:ALS-NCG} were
applied to the convex quadratic problem with an exact line
search and $P$ were represented by a preconditioning matrix
$\mathbf{P}$, then Algorithm \ref{Algorithm:ALS-NCG} would
be equivalent to preconditioned CG with preconditioner
$\mathbf{P}$. In \cite{DeSterck:2015}, an overview and an
in-depth analysis are given of the different possible forms
for $\overline{\beta}_{k+1}$, however \eqn{AALS:PRBeta} performed best in our numerical experiments.
Note that the primary computational cost in Algorithm \ref{Algorithm:ALS-NCG} comes from computing both the ALS iteration, $P(\mathbf{x}_k)$, and $\alpha_k$, the step length parameter using a line search.  

\begin{algorithm}[tb]
 \SetAlgoLined
 \KwIn{$\mathbf{x}_0$}
 \KwOut{$\mathbf{x}_k$}
 $\overline{\mathbf{g}}_0 \leftarrow \mathbf{x}_0 -
 P(\mathbf{x}_0)$\;
 $\mathbf{p}_0 \leftarrow -\overline{\mathbf{g}}_0$\;
 $k \leftarrow 0$\;
 \While{$\mathbf{g}_k \neq 0$}{
  Compute $\alpha_k$\;
  $\mathbf{x}_{k+1} \leftarrow \mathbf{x}_k + \alpha_k\mathbf{p}_k$\;
  $\overline{\mathbf{g}}_{k+1} \leftarrow \mathbf{x}_{k+1} - P(\mathbf{x}_{k+1})$\;
  Compute $\overline{\beta}_{k+1}$\;
  $\mathbf{p}_{k+1} \leftarrow -\overline{\mathbf{g}}_{k+1}+\overline{\beta}_{k+1}\mathbf{p}_k$\;
  $k \leftarrow k +1$\;
   }
 \caption{Accelerated ALS (ALS-NCG)}
 \label{Algorithm:ALS-NCG}
\end{algorithm}

\section{Serial Performance of ALS-NCG}\label{Section:NR}

The ALS and ALS-NCG algorithms were implemented in serial in MATLAB, and evaluated using the MovieLens 20M dataset \cite{MovieLens:2015}.
The entire MovieLens 20M dataset has 138,493 users, 27,278
movies,
and just over 20 million ratings, where each user has rated at least 20 movies.
To investigate the algorithmic performance as a function of problem size, both the ALS and ALS-NCG algorithms were run on subsets of the MovieLens 20M dataset. 
In creating subsets, we excluded outlier users with either very many or very few movie ratings relative to the median number of ratings per user. 
To construct a subset with $n_u$ users, the users from the full dataset were sorted in descending order by the values of $n_{u_i}$ for each user.
Denoting the index of the user with the median number of ratings by $c$, the set of users from index $c - \floor{\frac{n_u}{2}}$ to $c + (\ceil{\frac{n_u}{2}} - 1)$ were included.
Once the users were determined, the same process was used to select the movies, where the ratings per movie were computed only for the chosen users.

All serial experiments were performed on a linux workstation
with a quad-core 3.16 GHz processor (Xeon X5460) and 8 GB of
RAM. For the ALS-NCG algorithm, the Mor\'{e}-Thuente line
search algorithm from the Poblano toolbox \cite{Poblano} was
used to compute $\alpha_k$.
The line search parameters were as follows: $10^{-4}$ for the sufficient decrease condition tolerance, $10^{-2}$ for the curvature condition tolerance, and an initial step length of 1 and a maximum of 20 iterations.
The stopping criteria for both ALS and ALS-NCG were the maximum number of iterations as well as a desired tolerance value in the gradient norm normalized by the number of variables, $\frac{1}{N}\|\gvec_k\|$ for $N = n_f\times(n_u+n_m)$.
For both algorithms, a normalized gradient norm of less than
$10^{-6}$ was required within at most $10^4$ iterations.
In addition, ALS-NCG had a maximum number of allowed function evaluations in the line search equal to $10^7$.

The serial tests were performed on ratings matrices of 4 different sizes:  $n_u \times n_m = 400 \times 80$, $800 \times 160$, $1600 \times 320$ and $3200 \times 640$, where $n_u$ is the number of users and $n_m$ is the number of movies. 
For each ratings matrix, ALS and ALS-NCG were used to solve the optimization problem in (\ref{Equation:OptProblem}) with $\lambda = 0.1$ and $n_f = 10$.
The algorithms were each run using 20 different random starting iterates (which were the same for both algorithms) until one of the stopping criteria was reached.
Table \ref{Table:TimingResults} summarizes the timing results for different problem sizes, where the given times are written in the form $a \pm b$ where $a$ is the mean time in seconds and $b$ is the standard deviation about the mean.
Since computing the gradient is not explicitly required in the ALS algorithm, the computation time for the gradient norm was excluded from the timing results for the ALS algorithm.
Runs that did not converge based on the gradient norm tolerance were not included in the mean and standard deviation calculations, however the only run that did not converge to $\frac{1}{N}\|\gvec_k\| < 10^{-6}$ before reaching the maximum number of iterations was a single $1600\times320$ ALS run.
The large standard deviations in the timing measurements stem from the variation in the number of iterations required to reduce the gradient norm.
The fourth column of Table \ref{Table:TimingResults} shows the acceleration factor of ALS-NCG, computed as the mean time for convergence of ALS divided by the mean time for convergence of ALS-NCG. We see from this table that ALS-NCG significantly accelerates the ALS algorithm for all problem sizes. Similarly, Table \ref{Table:IterationResults} summarizes the number of iterations required to reach convergence for each algorithm. Again, the results were calculated based on converged runs only.

{\renewcommand{\arraystretch}{1.05}{
\begin{table}[b]
\centering
\caption{Timing results of ALS and ALS-NCG.}
\label{Table:TimingResults}
\begin{tabular}{c | c | c | c}
\hline
{Problem Size} & \multicolumn{2}{|c|}{Time (s)} & Acceleration\\
\cline{2-3}
{$n_u \times n_m$} & {ALS} & {ALS-NCG} & Factor\\
\hline
\hline
$400 \times 80$   & 56.50 $\pm$ 38.06  & 12.22 $\pm$  4.25 & 4.62\\
$800 \times 160$  & 162.0 $\pm$ 89.94  & 47.57 $\pm$ 20.02 & 3.41\\
$1600 \times 320$ & 30.61 $\pm$ 120.3  & 116.2 $\pm$ 31.56 & 2.84\\
$3200 \times 640$ & 960.8 $\pm$ 364.0  & 303.7 $\pm$ 111.3 & 3.16\\
\hline
\end{tabular}
\end{table}}

\begin{table}[b]
\centering
\caption{Iteration results of ALS and ALS-NCG.}
\label{Table:IterationResults}
\begin{tabular}{c | c | c | c}
\hline
{Problem Size} & \multicolumn{2}{|c|}{Number of Iterations} & Acceleration\\
\cline{2-3}
$n_u \times n_m$ & {ALS} & {ALS-NCG} & Factor \\
\hline
\hline
$400 \times 80$   & 2181 $\pm$ 1466 & 158.8 $\pm$  54.3 & 12.74\\
$800 \times 160$  & 3048 $\pm$ 1689 & 290.4 $\pm$ 128.1 & 10.50\\
$1600 \times 320$ & 3014 $\pm$ 1098 & 302.9 $\pm$  86.3 & 9.95\\
$3200 \times 640$ & 4231 $\pm$ 1602 & 329.6 $\pm$ 127.5 & 12.84\\
\hline
\end{tabular}
\end{table}

From Tables \ref{Table:TimingResults} and \ref{Table:IterationResults} it is clear that ALS-NCG accelerates the convergence of the ALS algorithm, using the gradient norm as the measure of convergence.  However, since the factor matrices $\mathbf{U}$ and $\mathbf{M}$ are used to make recommendations, we would also like to examine the convergence of the algorithms in terms of the accuracy of the resultant recommendations.  
In particular, we are interested in the rankings of the top $t$ movies  (e.g. top 20 movies) for each user, and want to explore how these rankings change with increasing number of iterations for both ALS and ALS-NCG. If the rankings of the top $t$ movies for each user no longer change, then the algorithm has likely computed an accurate solution.

To measure the relative difference in rankings we use a
metric that is based on the number of pairwise swaps required to
convert a vector of movie rankings $\mathbf{p}_2$ into another
vector $\mathbf{p}_1$, but only for the top $t$ movies. 
We use a modified Kendall-Tau \cite{Kendall:1938} distance to
compute the difference between ranking vectors based only on
the rankings of the top $t$ items, normalize the distance to range
in $[0,1]$, and subsequently average the distances over all users.
To illustrate the ranking metric,
consider the following example, where $\mathbf{p}_1 = 
[6,3,1,2,4,5]$, and  $\mathbf{p}_2 = 
[3,4,2,5,6,1]$ are two different rankings of movies for user $i$, and let $t = 2$.
We begin by finding the top $t$ movies in $\mathbf{p}_1$.
Here, user $i$ ranks movie 6 highest.
In $\mathbf{p}_2$, movie 6 is ranked 5th, and there are 4 pairwise inversions required to place movie 6 in the first component of $\pvec_2$, producing a new ranking vector for the second iteration, $\tilde{\pvec}_2 = [6,3,4,2,5,1]$.
The 2nd highest ranked movie in $\mathbf{p}_1$ is movie 3. 
In $\tilde{\pvec}_2$, movie 3 is already ranked 2nd; as such, no further inversions are required.
In this case, the total number of pairwise swaps required to match ${\mathbf{p}}_2$ to ${\mathbf{p}}_1$ for the top 2 rankings is $s_i = 4$.
The distance $s_i$ is normalized to 1 if no inversions were needed (i.e.~$\mathbf{p}_1 = \mathbf{p}_2$ in the top $t$ spots) and 0 if the maximum number of inversions are needed.
The maximum number of inversions occurs if $\mathbf{p}_2$ and $\mathbf{p}_1$ are in opposite order, requiring $n_m-1$ inversions in the first step, $n_m-2$ inversions in the second step, and so on until $n_m-t$ inversions are needed in the $t$-th step.
Thus, the maximum number of inversions is $s_{\text{max}} = (n_m-1) + (n_m-2) + \ldots + (n_m-t) = \frac{t}{2}(2n_m-t-1)$, yielding the ranking accuracy metric for user $i$ as 
\begin{equation*}
	q_i = 1 - \frac{s_i}{s_{\text{max}}}.
\end{equation*}
The total ranking accuracy for an algorithm is taken as the average value of $q_i$ across all users relative to ranking of the top $t$ movies for each user produced from the solution obtained after the algorithm converged.

\begin{table}[b]
\centering
\caption{Ranking accuracy timing results for problem size $400 \times 80$}
\label{Table:TimeAccuracy400x80}
\begin{tabular}{c | c | c | c}
\hline
{Ranking} & \multicolumn{2}{|c|}{Time (s)} & Acceleration\\
\cline{2-3}
Accuracy & ALS & ALS-NCG & Factor \\
\hline
\hline
$70\%$ &   0.37 $\pm$  0.17 &  0.65 $\pm$ 0.16 & 0.58\\
$80\%$ &   7.88 $\pm$  7.03 &  2.87 $\pm$ 1.72 & 2.74\\
$90\%$ &  37.08 $\pm$ 37.62 &  6.92 $\pm$ 4.47 & 5.36\\
$100\%$& 101.29 $\pm$ 68.04 & 14.07 $\pm$ 5.25 & 7.20\\
\hline
\end{tabular}
\end{table}

The normalized Kendall-Tau ranking metric for the top $t$ rankings described above was used to evaluate the accuracy of ALS and ALS-NCG as a function of running time for $t = 20$.
Tables \ref{Table:TimeAccuracy400x80} and \ref{Table:TimeAccuracy800x160} summarize the time needed for each algorithm to reach a specified percentage ranking accuracy for ratings matrices with different sizes.
Both algorithms were run with 20 different random starting iterates, and the time to reach a given ranking accuracy is written in the form $a \pm b$ where $a$ is the mean time across all converged runs and $b$ is the standard deviation about the mean.
In the fourth column of Tables \ref{Table:TimeAccuracy400x80} and \ref{Table:TimeAccuracy800x160}, we have computed the acceleration factor of ALS-NCG as the ratio of the mean time for ALS to the mean time for ALS-NCG. 
To reach $70\%$ accuracy, ALS is faster for both problem sizes, however for accuracies greater than 70\%, ALS-NCG converges in significantly less time than ALS.
Thus, if an accurate solution is desired, ALS-NCG is much faster in general than ALS.
This is further illustrated in \fig{AC}, which shows the average time needed for both algorithms to reach a given ranking accuracy for the $400\times80$ ratings matrix. 
Here, both ALS and ALS-NCG reach a ranking accuracy of $75\%$ in less than a second, however it then takes ALS approximately 100 s to increase the ranking accuracy from $75\%$ to $100\%$, while ALS-NCG reaches the final ranking in only 14 s.

\begin{figure}[t]
\centering
\epsfig{file=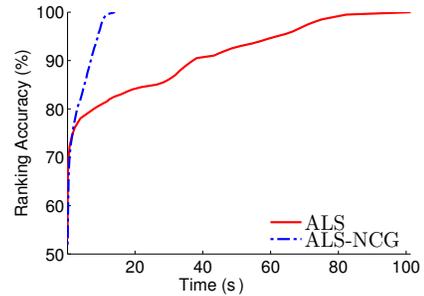,scale=0.4}
\caption{Average ranking accuracy versus time for problem size $400 \times 80$.}
\label{Figure:AC}
\end{figure}

\begin{table}[b]
\centering
\caption{Ranking accuracy timing results for problem size $800 \times 160$.}
\label{Table:TimeAccuracy800x160}
\begin{tabular}{c | c | c | c}
\hline
{Ranking} & \multicolumn{2}{|c|}{Time (s)} & Acceleration\\
\cline{2-3}
Accuracy & ALS & ALS-NCG & Factor\\
\hline
\hline
$70\%$ &   0.76 $\pm$   0.34  &  1.43 $\pm$  0.35 & 0.53\\
$80\%$ &  19.83 $\pm$  13.67  & 12.61 $\pm$ 10.13 & 1.57\\
$90\%$ & 103.91 $\pm$  70.97  & 33.25 $\pm$ 22.34 & 3.12\\
$100\%$& 310.98 $\pm$ 168.67  & 59.88 $\pm$ 28.20 & 5.19\\
\hline
\end{tabular}
\end{table}

\section{Parallel Implementation of ALS-NCG in Spark}\label{Section:Parallel}

\subsection{Apache Spark}
Apache Spark is a fault-tolerant, in-memory cluster computing framework designed to supersede MapReduce by maintaining program data in memory as much as possible between distributed operations.
The Spark environment is built upon two components: a data abstraction, termed a resilient distributed dataset (RDD) \cite{zaharia2012resilient}, and the task scheduler, which uses a \emph{delay scheduling} algorithm \cite{zaharia2010delay}.
We describe the fundamental aspects of RDDs and the scheduler below.

\textbf{Resilient Distributed Datasets}.
RDDs are immutable, distributed datasets that are evaluated lazily via their provenance information---that is, their functional relation to other RDDs or datasets in stable storage. 
To describe an RDD, consider an immutable distributed dataset $D$ of $k$ records with homogeneous type: $D = \bigcup_i^k d_i$ with $d_i \in \mathcal{D}$.
The distribution of $D$ across a computer network of nodes $\Set{v_{\alpha}}$, such that $d_i$ is stored in memory or on disk on node $v_{\alpha}$, is termed its \emph{partitioning} according to a partition function $P(d_i) = v_{\alpha}$.
If $D$ is expressible as a finite sequence of deterministic operations on other datasets $D_1,\ldots, D_l$ that are either RDDs or persistent records, then its lineage may be written as a directed acyclic graph $\lineage$ formed with the parent datasets $\Set{D_i}$ as the vertices, and the operations along the edges.
Thus, an RDD of type $\mathcal{D}$ (written $\rdd{\mathcal{D}}$) is the tuple $(D,P,\lineage)$.

Physically computing the records $\Set{d_i}$ of an RDD is termed its \emph{materialization}, and is managed by the Spark scheduler program.
To allocate computational tasks to the compute nodes, the scheduler traverses an RDD's lineage graph $\lineage$ and divides the required operations into stages of local computations on parent RDD partitions.
Suppose that $R_0 = (\bigcup_i x_i,P_0,\lineage_0)$ were an RDD of numeric type $\rdd{\mathbb{R}}$, and let $R_1 = (\bigcup_i y_i, P_1, \lineage_1)$ be the RDD resulting from the application of function $f:\mathbb{R}\rightarrow\mathbb{R}$ to each record of $R_0$. 
To compute $\Set{y_i}$, $R_1$ has only a single parent in the graph $\lineage_1$, and hence the set of tasks to perform is $\Set{f(x_i)}$.
This type of operation is termed a \emph{map} operation.
If $P_1 = P_0$, $\lineage_1$ is said to have a \emph{narrow} dependency on $R_0$: each $y_i$ may be computed locally from $x_i$, and the scheduler would allocate the task $f(x_i)$ to a node that stores $x_i$.

Stages consist only of local map operations, and are bounded by \emph{shuffle} operations that require communication and data transfer between the compute nodes.
For example, shuffling is necessary to perform \emph{reduce} operations on RDDs, wherein a scalar value is produced from an associative binary operator applied to each element of the dataset.  
In implementation, a shuffle is performed by writing the results of the tasks in the preceding stage, $\Set{f(x_i)}$, to a local file buffer. 
These shuffle files may or may not be written to disk, depending on the operating system's page table, and are fetched by remote nodes as needed in the subsequent stage.

\textbf{Delay Scheduling and Fault Tolerance}. 
The simple delay scheduling algorithm \cite{zaharia2010delay}
prioritizes data locality when submitting tasks to the available compute nodes.
If $v_{\alpha}$ stores the needed parent partition $x_i$ to compute task $f(x_i)$, but is temporarily unavailable due to faults or stochastic delays, rather than submitting the task on another node, the scheduler will wait until $v_{\alpha}$ is free.
However, if $v_{\alpha}$ does not become available within a specified maximum delay time (several seconds in practice), the scheduler will resubmit the tasks to a different compute node.
However, as $x_i$ is not available in memory on the different node, the lineage $\lineage_0$ of the RDD $R_0$ must be traversed further, and the tasks required to compute $x_i$ from the parent RDDs of $R_0$ will be submitted for computation in addition to $f(x_i)$. 
Thus, fault tolerance is achieved in the system through recomputation.

\subsection{ALS Implementation in Spark}
The Apache Spark codebase contains a parallel implementation of ALS for the collaborative filtering model of \cite{Zhou:2008,Koren:2009}; see \cite{Johnson2014}.
We briefly outline the main execution and relevant optimizations of the existing ALS implementation.
The implementation stores the factor matrices $\Umat$ and $\Mmat$ as single precision column block matrices, with each block as an RDD partition. A $\Umat$ column block stores factor vectors for a subset of users, and an $\Mmat$ column block stores factor vectors for a subset of movies.
The ratings matrix is stored twice: both $\Rmat$ and $\Rmat^T$ are stored in separate RDDs, partitioned in row blocks (i.e., $\Rmat$ is partitioned by users and $\Rmat^T$ by movies).
Row blocks of $\Rmat$ and $\Rmat^T$ are stored in a compressed sparse matrix format.
The $\Umat$ and $\Rmat$ RDDs have the same partitioning, such that a node that stores a partition of $\Umat$ also stores the partition of $\Rmat$ such that the ratings for each user in its partitions of $\Umat$ are locally available.
When updating a $\Umat$ block according to \eqn{SolveU}, the required ratings are available in a local $\Rmat$ block, but the movie factor vectors in $\Mmat$ corresponding to the movies rated by the users in a local $\Umat$ block must be shuffled across the network.
These movie factors are fetched from different nodes, and, as explained below, an optimized routing table strategy is used from \cite{Johnson2014} that avoids sending duplicate information. 
Similarly, updating a block of $\Mmat$ according to \eqn{SolveM} uses ratings data stored in a local $\Rmat^T$ block, but requires shuffling of $\Umat$ factor vectors using a second routing table.

\textbf{Block Partitioning}.
All RDDs are partitioned into $n_b$ partitions,
where $n_b$ is an integer multiple of the number of
available compute cores in practice.
For example, $\Mmat$ is divided into column blocks ${\Mmat_{j_b}}$ with block (movie) index $j_b \in \Set{0, \ldots, n_b-1}$ by hash partitioning the movie factor vectors such that $\mvec \in \Rmat_{j_b}$ if $j \equiv j_b \pmod{n_b}$ as in \fig{partition}.

	\begin{figure}[tb!]
		\centering
		\resizebox{1\columnwidth}{!}{
			\input{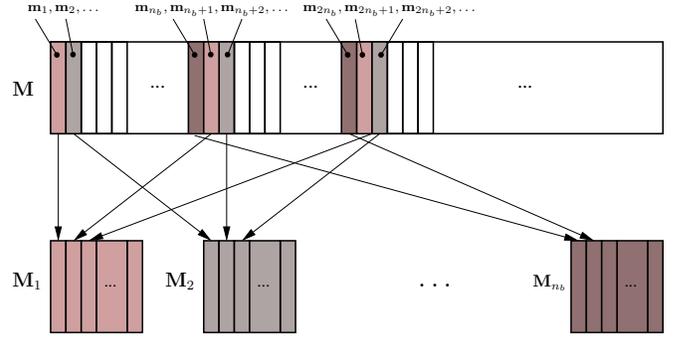}
		}
		\vspace{-3ex}
		\caption{
	Hash partitioning of the columns of $\Mmat$, depicted as rectangles, into $n_b$ blocks $\Mmat_1, \ldots, \Mmat_{n_b}$.
	Each block $M_{j_b}$ and its index $j_b$ forms a partition of the $\Mmat$ RDD of type $\rdd{(j_b,M_{j_b})}$.
}
		\label{Figure:partition}
	\end{figure}

Similarly, $\Umat$ is hash partitioned into column blocks ${\Umat_{i_b}}$ with block (user) index $i_b \in \Set{0, \ldots, n_b-1}$.
The RDDs for $\Mmat$ and $\Umat$ can be taken as type $\rdd{(j_b,\Mmat_{j_b})}$ and $\rdd{(i_b,\Umat_{i_b}^T)}$, where 
the blocks are tracked by the indices $j_b$ and $i_b$.
$\Rmat$ is partitioned by rows (users) into blocks with type $\rdd{(i_b,\Rmat_{i_b})}$ with the same partitioning as the RDD representing $\Umat$ (and similarly for the $\Rmat^T$ and $\Mmat$ RDDs).
By sharing the same user-based partitioning scheme, the blocks $\Rmat_{i_b}$ and $\Umat_{i_b}$ are normally located on the same compute node, except when faults occur. 
The same applies to $\Rmat^T_{j_b}$ and $\Mmat_{j_b}$ due to the movie-based partitioning scheme.

\textbf{Routing Table}.
\fig{shuffle} shows how a routing table optimizes data shuffling in ALS.
Suppose we want to update the user factor block $\Umat_{i_b}$ according to \eqn{SolveU}.
The required ratings data, $\Rmat_{i_b}$, is stored locally,
but a global shuffle is required to obtain all movie factor vectors in $\Mmat$ that correspond to the movies rated by the users in $\Umat_{i_b}$.
To optimize the data transfer, a routing table $T_m(\mvec)$ is constructed by determining, for each of the movie factor blocks $\Mmat_{j_b}$, which factor vectors have to be sent to each partition of $\Rmat_{i_b}$ (that may reside on other nodes).
In \fig{shuffle} (a), the blocks $\Mmat_{j_b}$ are \emph{filtered} using $T_m(\mvec)$ such that a given $\mvec \in \Mmat_{j_b}$ is written to the buffer destined for $\Rmat_{i_b}$, $\Mmat_{j_b}^{[i_b]}$, \emph{only once} regardless of how many $\uvec$ in $\Umat_{i_b}$ have ratings for movie $j$, and \emph{only} if there is at least one $\uvec$ in $\Umat_{i_b}$ that has rated movie $j$.
	This is shown by the hatching of each $\mvec$ vector in \fig{shuffle} (a); for instance, the first column in $\Mmat_1$ is written only to $\Mmat_1^{[1]}$ and has one set of hatching lines, but the last column is written to both $\Mmat_1^{[2]}$ and $\Mmat_{1}^{[n_b]}$ and correspondingly has two sets of hatching lines.
	Once the buffers are constructed, they are shuffled to the partitions of $\Rmat$, as in \fig{shuffle} (b) such that both the movie factors and ratings are locally available to compute the new $\Umat_{i_b}$ block, as in \fig{shuffle} (c).
The routing table formulation for shuffling $\Umat$, with mapping $T_u(\uvec)$, is analogous.
Note that the routing tables are constructed before the while loop in Algorithm \ref{Algorithm:ALS}, and hence do not need recomputation in each iteration.

	\begin{figure}[tb!]
		\centering
		\resizebox{1\columnwidth}{!}{
			\input{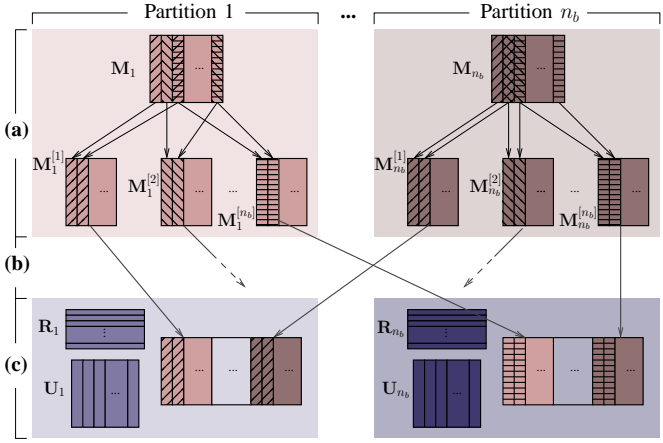}
		}
		\vspace{-3ex}
		\caption{
	Schematic of the use of the routing table $T_m(m_j)$ in the Spark ALS shuffle.
	In (a) the blocks $\Set{\Mmat_{j_b}}$ are filtered using $T_m(\mvec)$ for each destination $\Rmat_{i_b}$ and shuffled to the respective blocks in (b), where arrows between the shaded backgrounds represent network data transfer between different compute nodes.
		In (c), when updating block $\Umat_{i_b}$, the ratings information is locally available in $\Rmat_{i_b}$.
}
		\label{Figure:shuffle}
	\end{figure}

\textbf{Evaluation of $\uvec$ and $\mvec$ via \eqn{SolveU} and \eqn{SolveM}}.
As in \fig{shuffle} (c), a compute node that stores $\Rmat_{i_b}$, 
will obtain $n_b$ buffered arrays of filtered movie factors. 
Once the factors have been shuffled, $\Amat_j$ is computed for each $\uvec$ as $\sum_{j\in\Iset_i}\mvec\mvec^T + n_{u_i}\Imat$, and $\vvec_j$ as $\sum_{j \in \Iset_i} r_{ij} \mvec$ using the Basic Linear Algebra Subprograms (BLAS) library \cite{blackford2002updated}. 
The resulting linear system for $\uvec$ is then solved via the Cholesky decomposition using LAPACK routines \cite{anderson1999lapack}, giving the computation to update $\Umat$ an asymptotic complexity of $\order{n_u n_f^3}$ since $n_u$ linear systems must be solved.
Solving for $\mvec$ is an identical operation with the appropriate routing table and has $\order{n_m n_f^3}$ complexity.

\subsection{ALS-NCG Implementation in Spark}
We now discuss our contributions in parallelizing ALS-NCG for Spark. 
Since the calculation of $\alpha_k$ in Algorithm
\ref{Algorithm:ALS-NCG} requires a line search, the main
technical challenges we address are how to compute the loss
function in \eqn{ActualLoss} and its gradient in an efficient way in Spark, obtaining good parallel performance.
To this end, we formulate a backtracking line search procedure that dramatically reduces the cost of multiple function evaluations, and we take advantage of the routing table mechanism to obtain fast communication.
We also extend the ALS implementation in Spark to support the additional NCG vector operations required in Algorithm \ref{Algorithm:ALS-NCG} using BLAS routines.

\textbf{Vector Storage}.
The additional vectors $\xpvec, \gvec, \gpvec$, and $\pvec$
were each split into two separate RDDs, such that blocks
corresponding to the components of $\uvec$ (see \eqn{x})
were stored in one RDD and partitioned in the same way as $\Umat$ with block index $i_b$.
Analogously, blocks corresponding to components of $\mvec$ were stored in another RDD, partitioned in the same way as $\Mmat$ with block index $j_b$.
This ensured that all vector blocks were aligned component-wise for vector operations; furthermore, the $\pvec$ blocks could also be shuffled efficiently using the routing tables in the line search (see below).

\textbf{Vector Operations}.
RDDs have a standard operation termed a \emph{join}, in which two RDDs $R_1$ and $R_2$ representing the datasets of tuples $\bigcup_i \left(k_i,a_i\right)$ and $\bigcup_i (k_i,b_i)$ with type $(\key,\A)$ and $(\key,\B)$, respectively, are combined to produce an RDD $R_3$ of type $\rdd{\left(\key,(\A,\B)\right)}$, where $R_3 = R_1.\text{join}(R_2)$ represents the dataset $\bigcup_l (k_l, (a_l,b_l))$ of combined tuples and $k_l$ is a key common to both $R_1$ and $R_2$. 
Parallel vector operations between RDD representations of blocked vectors $\xvec$ and $\yvec$ were implemented by joining the RDDs by their block index $i_b$ and calling BLAS level 1 interfaces on the vectors within the resultant tuples of aligned vector blocks, $\Set{(\xvec_{i_b},\yvec_{i_b})}$. 
Since the RDD implementations of vectors had the same partitioning schemes, this operation was local to a compute node, and hence very inexpensive.
One caveat, however, is that BLAS subprograms generally perform modifications to vectors \emph{in place}, overwriting their previous components.
For fault tolerance, RDDs must be immutable; as such, whenever BLAS operations were required on the records of an RDD, an entirely new vector was allocated in memory and the RDD's contents copied to it.
Algorithm \ref{Algorithm:rddAxpby} shows the operations required for the vector addition $a\xvec + b\yvec$ using the \texttt{BLAS.axpby} routine (note that the result overwrites $\yvec$ in place).
The inner product of two block vector RDDs and norm of a
single block vector RDD were implemented in a similar manner
to vector addition, with an additional reduce operation
to sum the scalar component-wise dot products.
These vector operations were used to
compute $\betap$ in \eqn{AALS:PRBeta}.

\begin{algorithm}[tb]
	{
		\small
 \SetAlgoNoLine
 \SetKwInOut{Input}{Input}
 \SetKwInOut{Output}{Output}
 \KwIn{$\xvec = \rdd{(i_b,\xvec_{i_b})}$; $\yvec = \rdd{(i_b,\yvec_{i_b})}$; $a,b \in \mathbb{R}$}
 \KwOut{$\zvec = a\xvec + b\yvec$}
 	$\zvec \leftarrow \xvec$.join($\yvec$).map\{ \\
	\Indp

	Allocate $\zvec_{i_b};$ \\
	$\zvec_{i_b} \leftarrow \yvec_{i_b};$ \\
	Call \texttt{BLAS.axpby}($a,\xvec_{i_b}, b,\zvec_{i_b}$);\\
	Yield $(i_b,\zvec_{i_b});$

 \Indm
 }

 \}
 \caption{RDD Block Vector \texttt{BLAS.axpby}}
 \label{Algorithm:rddAxpby}
\end{algorithm}

\textbf{Line Search \& Loss Function Evaluation}.
We present a computationally cheap way to implement a
backtracking line search in Spark for minimizing \eqn{ActualLoss}.
A backtracking line search minimizes a function
$f(\xvec)$ along a descent direction $\pvec$ by decreasing
the step size in each iteration by a factor of $\tau \in (
0,1 )$, and terminates once a sufficient decrease in $f(\xvec)$ is achieved.
This simple procedure is summarized in Algorithm
\ref{Algorithm:ls}, where it is important to note that a
line search requires computing $\gvec$ once, as well as $f(\xvec)$ in each iteration
of the line search. 
To avoid multiple shuffles in each line search iteration,
instead of performing multiple evaluations of
\eqn{ActualLoss} directly, we
constructed a polynomial with degree 4 in the step size $\alpha$
by expanding \eqn{ActualLoss} with $\uvec = \xvec_{u_i} +
\alpha \pvec_{u_i}$ and $\mvec = \xvec_{m_j} + \alpha
\pvec_{m_j}$, where $\xvec_{u_i}$ refers to the components of $\xvec$ related to user $i$ and mutatis mutandis for the vectors $\pvec_{u_i}$, $\xvec_{m_j}$, and $\pvec_{m_j}$.
From the bilinearity of the inner product,
this polynomial has the form
$	Q(\alpha) = \sum_{n = 0}^4 \left(\sum_{(i,j) \in \mathcal{I}} C^{[n]}_{ij}\right) \alpha^n$,
where the terms $C^{[n]}_{ij}$ in the summation for each coefficient 
only require level 1 BLAS operations between the block
vectors $\xvec_{u_i},\pvec_{u_i},\xvec_{m_j}$, and
$\pvec_{m_j}$ for known $(i,j)\in\Iset$.
Thus, coefficients of $Q(\alpha)$ were computed at the
beginning of the line search with a
\emph{single}
shuffle operation using
the routing table $T_m(\mvec)$ 
to match vector pairs with dot products contributing to the
coefficients. 
Here, $T_m(\mvec)$ was chosen since there is far less communication required to
shuffle $\Set{\mvec}$, as $n_u \gg n_m$.
Since each iteration of the line search was very fast after computing the coefficients of
$Q(\alpha)$,
we used relatively large values of $\tau = 0.9$, $c = 0.5$, and $\alpha_0 = 10$ in Algorithm \ref{Algorithm:ls} that searched intensively along direction $\pvec$.
\begin{algorithm}[b]
 \SetAlgoLined
 \KwIn{$\xvec$, $\pvec$, $\gvec$, $\alpha_0$, $c \in ( 0,1 )$, $\tau \in ( 0,1 )$}
 \KwOut{$\alpha_k$}

 $k \leftarrow 0;$

 \While{$f(\xvec + \alpha_k \pvec) - f(\xvec) > \alpha_k\,c\,\gvec^T \pvec $}{
	 $\alpha_{k+1} \leftarrow \tau\, \alpha_{k};$\\
	 $k \leftarrow k+1;$
 }
 \caption{Backtracking Line Search}
 \label{Algorithm:ls}
\end{algorithm}

\textbf{Gradient Evaluation}. 
We computed $\gvec_k$ with respect to a block for $\uvec$ using only BLAS level 1 operations as $2\lambda n_{u_i}\uvec + 2\sum_{j\in\Iset_i} \mvec (\uvec^T\mvec - r_{ij})$, with an analogous operation with respect to each $\mvec$.
As this computation requires matching the $\uvec$ and $\mvec$ factors, the routing tables $T_u(\uvec)$ and $T_m(\mvec)$ were used to shuffle $\uvec$ and $\mvec$ consecutively.
Evaluating $\gvec_k$ can be performed with $\order{n_f (n_{u}\sum_i n_{u_i} + n_m \sum_j n_{m_j})}$ operations, but requires two shuffles.
As such, the gradient computation required as much communication as a single iteration of ALS in which both $\Umat$ and $\Mmat$ are updated.

\section{ Parallel Performance of ALS-NCG}
\label{Section:sparkExperiments}

Our comparison tests of the ALS and ALS-NCG algorithms in Spark were performed on a computing cluster composed of 16 homogeneous compute nodes, 1 storage node hosting a network filesystem, and 1 head node.
The nodes were interconnected by a 10 Gb ethernet managed switch (PowerConnect 8164).
Each compute node was a 64 bit rack server (PowerEdge R620)
running Ubuntu 14.04, with linux kernel 3.13.
The compute nodes all had two 8-core 2.60 GHz chips (Xeon E5-2670) and 256 GB of SDRAM. 
The head node had the same processors as the compute
nodes, but had 512 GB of RAM.
The single storage node (PowerEdge R720) contained two 2 GHz processors, each with 6 cores (Xeon E5-2620), 64 GB of memory, and 12 hard disk drives of 4 TB capacity and 7200 RPM nominal speed.
Finally, compute nodes were equipped with 6 ext4-formatted local SCSI 10k RPM hard disk drives, each with a 600 GB capacity.

Our Apache Spark assembly was built from a snapshot of the 1.3 release using Oracle's Java 7 distribution, Scala 2.10, and Hadoop 1.0.4.
Input files to Spark programs were stored on the storage node in plain text.
The SCSI hard drives on the compute nodes' local filesystems were used as Spark spilling and scratch directories, and the Spark checkpoint directory for persisting RDDs was specified in the network filesystem hosted by the storage node, and accessible to all nodes in the cluster.
Shuffle files were consolidated into larger files, as recommended for ext4 filesystems \cite{davidson2013optimizing}.
In our experiments, the Spark master was executed on the head node, and a single instance of a Spark executor was created on each compute node.
It was empirically found that the ideal number of cores to make available to the Spark driver was 16 per node, or the number of \emph{physical} cores for a total of 256 available cores.
The value of $n_b$ was set to the number of cores in all experiments.

To compare the performance of ALS and ALS-NCG in Spark, the two implementations were tested on the MovieLens 20M dataset with $\lambda = 0.01$ and $n_f = 100$.
In both algorithms, the RDDs were checkpointed to persistent storage every 10 iterations, since it is a widely known issue in the Spark community that RDDs with very long lineage graphs cause stack overflow errors when the scheduler recursively traverses their lineage \cite{Dai2015}.
Since RDDs are materialized via lazy evaluation, to obtain timing measurements, actions were triggered to physically compute the partitions of each block vector RDD at the end of each iteration in Algorithm \ref{Algorithm:ALS} and \ref{Algorithm:ALS-NCG}.
For each experimental run of ALS and ALS-NCG, two experiments with the same initial user and movie factors were performed: in one, the gradient norm in each iteration was computed and printed (incurring additional operations); in the other experiment, no additional computations were performed such that the elapsed times for each iteration were correctly measured.

\fig{convergence} shows the convergence in gradient norm, normalized by the degrees of freedom $N = n_f\times(n_u+n_m)$, for six separate runs of ALS and ALS-NCG on 8 compute nodes for the MovieLens 20M dataset.
The subplots (a) and (b) show $\frac{1}{N}\|\mathbf{g}_k\|$ over 100 iterations and 25 minutes, respectively.
This time frame was chosen since it took just over 20 minutes for ALS-NCG to complete 100 iterations; note that in \fig{convergence} (b), 200 iterations of ALS are shown, since with this problem size it took approximately twice as long to run a single iteration of ALS-NCG.
The shaded regions in \fig{convergence} show the standard deviation about the mean value for $\frac{1}{N}\|\gvec_k\|$ across all runs, computed for non-overlapping windows of 3 iterations for subplot (a) and 30 s for subplot (b).
Even within the uncertainty bounds, ALS-NCG requires much less time and many fewer iterations
than ALS to reach accurate values
of $\frac{1}{N}\norm{\gvec_k}$ (e.g. below $10^{-3}$).

	\begin{figure}[tb!]
		\centering
		\resizebox{0.80\columnwidth}{!}{
			{
				\footnotesize
				\input{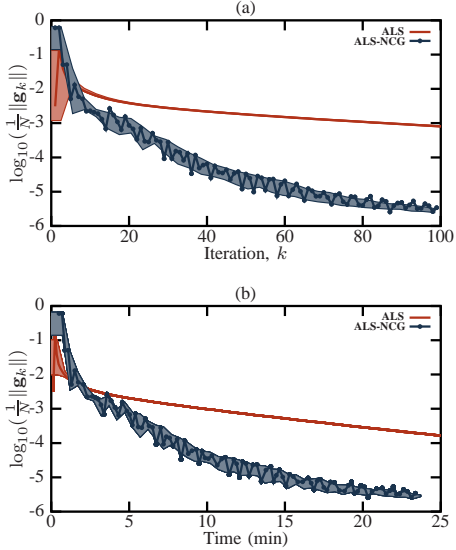}
			}
		}
		\vspace{-3ex}
		\caption{
	Convergence in normalized gradient norm $\frac{1}{N}\|\gvec_k\|$ for 6 instances of ALS and ALS-NCG with different starting values in both (a) iteration and (b) clock time, for the MovieLens 20M dataset. 
	The two solid lines in each panel show actual convergence traces for one of the instances, while the shaded regions show the standard deviation about the mean value over all instances, computed for non-overlapping windows of 3 iterations for (a), and 30 s for (b). 
	The experiments were conducted on 8 nodes (128 cores).
}
		\label{Figure:convergence}
	\end{figure}

The operations that we have implemented in ALS-NCG that are additional to standard ALS in Spark have computational complexity that is linear in problem size.
To verify the expected linear scaling, experiments with a constant number of nodes and increasing problem size up to 800 million ratings were performed on 16 compute nodes.
Synthetic ratings matrices were constructed by sampling the MovieLens 20M dataset such that the synthetic dataset had the same statistical sparsity, realistically simulating the data transfer patterns in each iteration.
To do this, first the dimension $n_u$ of the sampled dataset was fixed, and for each user $n_{u_i}$ was sampled from the empirical probability distribution $p(n_{u_i}|\Rmat)$ of how many movies each user ranked, computed empirically from the MovieLens 20M dataset.
The $n_{u_i}$ movies were then sampled from the empirical likelihood of sampling the $j$th movie, $p(\mvec|\Rmat)$, and the resultant rating value was sampled from the distribution of numerical values for all ratings.
The choice of scaling up the users (for a fixed set of movies) was made
to model the situation in which an industry's user base grows far more rapidly than its items.

\fig{gran_scaling} shows the linear scaling in computation time for both ALS and ALS-NCG.
The values shown are average times per iteration over 50 iterations, for $n_u$ from 1 to 6 million, corresponding to the range from 133 to 800 million ratings.  
The error bars show the uncertainty in this measurement, where the standard deviation takes into account that
iterations with and without checkpointing come from two different populations with different average times.
The uncertainty in the time per iteration for ALS-NCG is larger due to the greater overhead of memory management and garbage collection by the Java Virtual Machine required with more RDDs.
While each iteration of ALS-NCG takes longer due to the additional line search and gradient computations, we note that many fewer iterations are required to converge.

	\begin{figure}[tb!]
		\centering
		\resizebox{0.80\columnwidth}{!}{
			\input{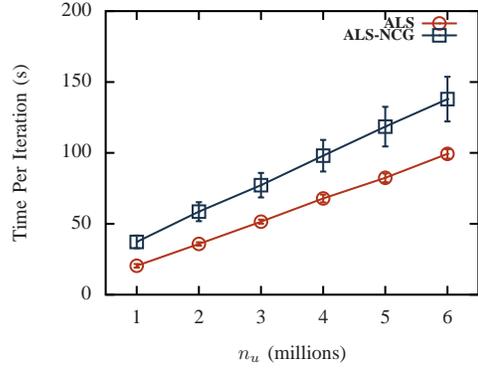}
		}
		\vspace{-3ex}
		\caption{
	Linear scaling of computation time per iteration with increasing $n_u$ on a synthetic dataset for ALS and ALS-NCG on 16 compute nodes (256 cores) for up to 6M users, corresponding to 800M ratings. 
}
		\label{Figure:gran_scaling}
	\end{figure}

Finally, we compute the relative speedup that was attainable on the large synthetic datasets.
For the value of $\frac{1}{N}\|\gvec_k\|$ in each iteration of ALS-NCG, we determined how many iterations of regular ALS were required to achieve an equal or lesser value gradient norm.  Due to the local variation in $\frac{1}{N}\|\gvec_k\|$ (as in \fig{convergence}), a moving average filter over every two iterations was applied to the ALS-NCG gradient norm values.
The total time required for ALS and ALS-NCG to reach a given gradient norm was then estimated from the average times per iteration in \fig{gran_scaling}. 
The ratios of these total times for ALS and ALS-NCG are shown in \fig{speedup} as the relative speedup factor for the 1M, 3M, and 6M users ratings matrices.
When an accurate solution is desired, ALS-NCG often achieves
faster convergence by a factor of 3 to 5, with the
acceleration factor increasing with greater desired accuracy
in the solution.

	\begin{figure}[tb!]
		\centering
		\resizebox{0.85\columnwidth}{!}{
			\input{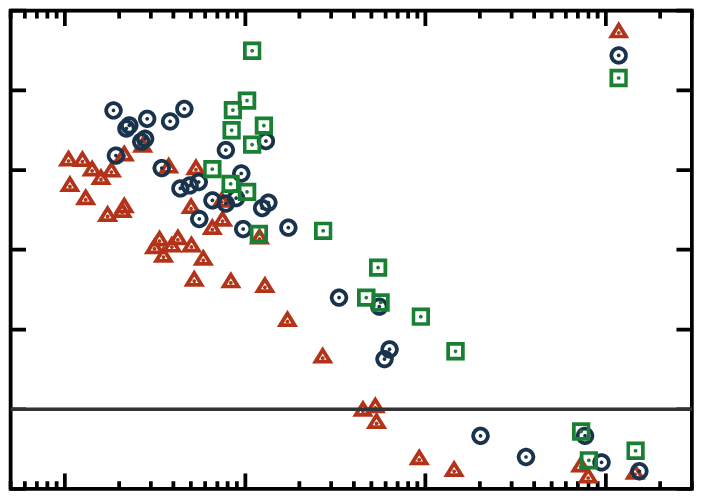}
		}
		\vspace{-3ex}
		\caption{
	Speedup of ALS-NCG over ALS as a function of normalized gradient norm on 16 compute nodes (256 cores), for a synthetic problem with up to 6M users and 800M ratings.
	ALS-NCG can easily outperform ALS by a factor of 4 or more, especially when accurate solutions (small normalized gradients) are required or problem sizes are large.
}
		\label{Figure:speedup}
	\end{figure}

\section{ Conclusion}
\label{Section:Conclusion}

In this paper, we have demonstrated how the addition of a nonlinear conjugate gradient wrapper can accelerate the convergence of ALS-based collaborative filtering algorithms when accurate solutions are desired.
Furthermore, our parallel ALS-NCG implementation can significantly speed up big data recommendation in the Apache Spark distributed computing environment, and the proposed NCG acceleration can be naturally extended to Hadoop or MPI environments. It may also be applicable to other optimization methods for collaborative filtering.
Though we have focused on a simple latent factor model, our acceleration can be used with any collaborative filtering model that uses ALS. 
We expect that our acceleration approach will be especially useful for advanced collaborative filtering models that achieve low root mean square error (RMSE), since these models require solving the optimization problem accurately, and that is precisely where accelerated ALS-NCG shows the most benefit over standalone ALS.
\section*{Acknowledgements}
This work was supported in part by 
NSERC of Canada, and by the Scalable Graph Factorization LDRD Project, 13-ERD-072, under the auspices of the U.S.~Department of Energy by Lawrence Livermore National Laboratory under Contract DE-AC52-07NA27344.

\insertBibliography

\end{document}

%% file: fig/partition.tex
\begin{picture}(0,0)%
\includegraphics{partition.eps}%
\end{picture}%
\setlength{\unitlength}{3947sp}%
\begingroup\makeatletter\ifx\SetFigFont\undefined%
\gdef\SetFigFont#1#2#3#4#5{%
  \reset@font\fontsize{#1}{#2pt}%
  \fontfamily{#3}\fontseries{#4}\fontshape{#5}%
  \selectfont}%
\fi\endgroup%
\begin{picture}(6412,3319)(3211,-2783)
\put(8176,-361){\makebox(0,0)[lb]{\smash{{\SetFigFont{12}{14.4}{\familydefault}{\mddefault}{\updefault}{\color[rgb]{0,0,0}...}%
}}}}
\put(4576,-361){\makebox(0,0)[lb]{\smash{{\SetFigFont{12}{14.4}{\familydefault}{\mddefault}{\updefault}{\color[rgb]{0,0,0}...}%
}}}}
\put(6076,-361){\makebox(0,0)[lb]{\smash{{\SetFigFont{12}{14.4}{\familydefault}{\mddefault}{\updefault}{\color[rgb]{0,0,0}...}%
}}}}
\put(4426,389){\makebox(0,0)[lb]{\smash{{\SetFigFont{8}{9.6}{\familydefault}{\mddefault}{\updefault}{\color[rgb]{0,0,0}$\mathbf{m}_{n_b},\mathbf{m}_{n_b + 1},\mathbf{m}_{n_b + 2}, \ldots$}%
}}}}
\put(6076,389){\makebox(0,0)[lb]{\smash{{\SetFigFont{8}{9.6}{\familydefault}{\mddefault}{\updefault}{\color[rgb]{0,0,0}$\mathbf{m}_{2n_b},\mathbf{m}_{2n_b + 1},\mathbf{m}_{2n_b + 2}, \ldots$}%
}}}}
\put(3376,389){\makebox(0,0)[lb]{\smash{{\SetFigFont{8}{9.6}{\familydefault}{\mddefault}{\updefault}{\color[rgb]{0,0,0}$\mathbf{m}_1,\mathbf{m}_2, \ldots$}%
}}}}
\put(7201,-2311){\makebox(0,0)[lb]{\smash{{\SetFigFont{17}{20.4}{\familydefault}{\mddefault}{\updefault}{\color[rgb]{0,0,0}$\ldots$}%
}}}}
\put(3226,-2311){\makebox(0,0)[lb]{\smash{{\SetFigFont{12}{14.4}{\familydefault}{\mddefault}{\updefault}{\color[rgb]{0,0,0}$\mathbf{M}_1$}%
}}}}
\put(9226,-2311){\makebox(0,0)[lb]{\smash{{\SetFigFont{10}{12.0}{\familydefault}{\mddefault}{\updefault}{\color[rgb]{0,0,0}...}%
}}}}
\put(8326,-2311){\makebox(0,0)[lb]{\smash{{\SetFigFont{10}{12.0}{\familydefault}{\mddefault}{\updefault}{\color[rgb]{0,0,0}$\mathbf{M}_{n_b}$}%
}}}}
\put(4726,-2311){\makebox(0,0)[lb]{\smash{{\SetFigFont{12}{14.4}{\familydefault}{\mddefault}{\updefault}{\color[rgb]{0,0,0}$\mathbf{M}_2$}%
}}}}
\put(5626,-2311){\makebox(0,0)[lb]{\smash{{\SetFigFont{10}{12.0}{\familydefault}{\mddefault}{\updefault}{\color[rgb]{0,0,0}...}%
}}}}
\put(4126,-2311){\makebox(0,0)[lb]{\smash{{\SetFigFont{10}{12.0}{\familydefault}{\mddefault}{\updefault}{\color[rgb]{0,0,0}...}%
}}}}
\put(3226,-436){\makebox(0,0)[lb]{\smash{{\SetFigFont{12}{14.4}{\familydefault}{\mddefault}{\updefault}{\color[rgb]{0,0,0}$\mathbf{M}$}%
}}}}
\end{picture}%

%% file: fig/shuffle.tex
\begin{picture}(0,0)%
\includegraphics{shuffle.eps}%
\end{picture}%
\setlength{\unitlength}{3947sp}%
\begingroup\makeatletter\ifx\SetFigFont\undefined%
\gdef\SetFigFont#1#2#3#4#5{%
  \reset@font\fontsize{#1}{#2pt}%
  \fontfamily{#3}\fontseries{#4}\fontshape{#5}%
  \selectfont}%
\fi\endgroup%
\begin{picture}(8802,5787)(1186,-7723)
\put(2476,-4411){\makebox(0,0)[lb]{\smash{{\SetFigFont{10}{12.0}{\familydefault}{\mddefault}{\updefault}{\color[rgb]{0,0,0}...}%
}}}}
\put(7051,-4411){\makebox(0,0)[lb]{\smash{{\SetFigFont{10}{12.0}{\familydefault}{\mddefault}{\updefault}{\color[rgb]{0,0,0}...}%
}}}}
\put(8326,-4411){\makebox(0,0)[lb]{\smash{{\SetFigFont{10}{12.0}{\familydefault}{\mddefault}{\updefault}{\color[rgb]{0,0,0}...}%
}}}}
\put(8326,-4411){\makebox(0,0)[lb]{\smash{{\SetFigFont{10}{12.0}{\familydefault}{\mddefault}{\updefault}{\color[rgb]{0,0,0}...}%
}}}}
\put(3751,-4411){\makebox(0,0)[lb]{\smash{{\SetFigFont{10}{12.0}{\familydefault}{\mddefault}{\updefault}{\color[rgb]{0,0,0}...}%
}}}}
\put(7126,-6361){\rotatebox{90.0}{\makebox(0,0)[lb]{\smash{{\SetFigFont{10}{12.0}{\familydefault}{\mddefault}{\updefault}{\color[rgb]{0,0,0}...}%
}}}}}
\put(7201,-7111){\makebox(0,0)[lb]{\smash{{\SetFigFont{10}{12.0}{\familydefault}{\mddefault}{\updefault}{\color[rgb]{0,0,0}...}%
}}}}
\put(2551,-6361){\rotatebox{90.0}{\makebox(0,0)[lb]{\smash{{\SetFigFont{10}{12.0}{\familydefault}{\mddefault}{\updefault}{\color[rgb]{0,0,0}...}%
}}}}}
\put(2626,-7111){\makebox(0,0)[lb]{\smash{{\SetFigFont{10}{12.0}{\familydefault}{\mddefault}{\updefault}{\color[rgb]{0,0,0}...}%
}}}}
\put(3751,-2761){\makebox(0,0)[lb]{\smash{{\SetFigFont{10}{12.0}{\familydefault}{\mddefault}{\updefault}{\color[rgb]{0,0,0}...}%
}}}}
\put(3751,-6811){\makebox(0,0)[lb]{\smash{{\SetFigFont{10}{12.0}{\familydefault}{\mddefault}{\updefault}{\color[rgb]{0,0,0}...}%
}}}}
\put(4201,-6811){\makebox(0,0)[lb]{\smash{{\SetFigFont{10}{12.0}{\familydefault}{\mddefault}{\updefault}{\color[rgb]{0,0,0}...}%
}}}}
\put(8326,-2761){\makebox(0,0)[lb]{\smash{{\SetFigFont{10}{12.0}{\familydefault}{\mddefault}{\updefault}{\color[rgb]{0,0,0}...}%
}}}}
\put(8326,-6811){\makebox(0,0)[lb]{\smash{{\SetFigFont{10}{12.0}{\familydefault}{\mddefault}{\updefault}{\color[rgb]{0,0,0}...}%
}}}}
\put(8776,-6811){\makebox(0,0)[lb]{\smash{{\SetFigFont{10}{12.0}{\familydefault}{\mddefault}{\updefault}{\color[rgb]{0,0,0}...}%
}}}}
\put(9526,-6811){\makebox(0,0)[lb]{\smash{{\SetFigFont{10}{12.0}{\familydefault}{\mddefault}{\updefault}{\color[rgb]{0,0,0}...}%
}}}}
\put(4951,-6811){\makebox(0,0)[lb]{\smash{{\SetFigFont{10}{12.0}{\familydefault}{\mddefault}{\updefault}{\color[rgb]{0,0,0}...}%
}}}}
\put(2626,-2836){\makebox(0,0)[lb]{\smash{{\SetFigFont{14}{16.8}{\familydefault}{\mddefault}{\updefault}{\color[rgb]{0,0,0}$\mathbf{M}_1$}%
}}}}
\put(7201,-2836){\makebox(0,0)[lb]{\smash{{\SetFigFont{14}{16.8}{\familydefault}{\mddefault}{\updefault}{\color[rgb]{0,0,0}$\mathbf{M}_{n_b}$}%
}}}}
\put(4201,-4411){\makebox(0,0)[lb]{\smash{{\SetFigFont{10}{12.0}{\familydefault}{\mddefault}{\updefault}{\color[rgb]{0,0,0}...}%
}}}}
\put(8776,-4411){\makebox(0,0)[lb]{\smash{{\SetFigFont{10}{12.0}{\familydefault}{\mddefault}{\updefault}{\color[rgb]{0,0,0}...}%
}}}}
\put(6151,-4111){\makebox(0,0)[lb]{\smash{{\SetFigFont{14}{16.8}{\familydefault}{\mddefault}{\updefault}{\color[rgb]{0,0,0}$\mathbf{M}_{n_b}^{[1]}$}%
}}}}
\put(4051,-4861){\makebox(0,0)[lb]{\smash{{\SetFigFont{14}{16.8}{\familydefault}{\mddefault}{\updefault}{\color[rgb]{0,0,0}$\mathbf{M}_1^{[n_b]}$}%
}}}}
\put(2851,-4411){\makebox(0,0)[lb]{\smash{{\SetFigFont{14}{16.8}{\familydefault}{\mddefault}{\updefault}{\color[rgb]{0,0,0}$\mathbf{M}_1^{[2]}$}%
}}}}
\put(7426,-4411){\makebox(0,0)[lb]{\smash{{\SetFigFont{14}{16.8}{\familydefault}{\mddefault}{\updefault}{\color[rgb]{0,0,0}$\mathbf{M}_{n_b}^{[2]}$}%
}}}}
\put(1576,-4111){\makebox(0,0)[lb]{\smash{{\SetFigFont{14}{16.8}{\familydefault}{\mddefault}{\updefault}{\color[rgb]{0,0,0}$\mathbf{M}_1^{[1]}$}%
}}}}
\put(8626,-4861){\makebox(0,0)[lb]{\smash{{\SetFigFont{14}{16.8}{\familydefault}{\mddefault}{\updefault}{\color[rgb]{0,0,0}$\mathbf{M}_{n_b}^{[n_b]}$}%
}}}}
\put(5701,-2086){\makebox(0,0)[lb]{\smash{{\SetFigFont{17}{20.4}{\rmdefault}{\bfdefault}{\updefault}{\color[rgb]{0,0,0}...}%
}}}}
\put(1651,-6286){\makebox(0,0)[lb]{\smash{{\SetFigFont{14}{16.8}{\familydefault}{\mddefault}{\updefault}{\color[rgb]{0,0,0}$\mathbf{R}_{1}$}%
}}}}
\put(1726,-7111){\makebox(0,0)[lb]{\smash{{\SetFigFont{14}{16.8}{\familydefault}{\mddefault}{\updefault}{\color[rgb]{0,0,0}$\mathbf{U}_{1}$}%
}}}}
\put(6191,-6295){\makebox(0,0)[lb]{\smash{{\SetFigFont{14}{16.8}{\familydefault}{\mddefault}{\updefault}{\color[rgb]{0,0,0}$\mathbf{R}_{n_b}$}%
}}}}
\put(6253,-7111){\makebox(0,0)[lb]{\smash{{\SetFigFont{14}{16.8}{\familydefault}{\mddefault}{\updefault}{\color[rgb]{0,0,0}$\mathbf{U}_{n_b}$}%
}}}}
\put(7576,-2086){\makebox(0,0)[lb]{\smash{{\SetFigFont{17}{20.4}{\familydefault}{\mddefault}{\updefault}{\color[rgb]{0,0,0}Partition $n_b$}%
}}}}
\put(3076,-2086){\makebox(0,0)[lb]{\smash{{\SetFigFont{17}{20.4}{\familydefault}{\mddefault}{\updefault}{\color[rgb]{0,0,0}Partition 1}%
}}}}
\put(1201,-3661){\makebox(0,0)[lb]{\smash{{\SetFigFont{17}{20.4}{\rmdefault}{\bfdefault}{\updefault}{\color[rgb]{0,0,0}(a)}%
}}}}
\put(1201,-5461){\makebox(0,0)[lb]{\smash{{\SetFigFont{17}{20.4}{\rmdefault}{\bfdefault}{\updefault}{\color[rgb]{0,0,0}(b)}%
}}}}
\put(1201,-6811){\makebox(0,0)[lb]{\smash{{\SetFigFont{17}{20.4}{\rmdefault}{\bfdefault}{\updefault}{\color[rgb]{0,0,0}(c)}%
}}}}
\put(5026,-4411){\makebox(0,0)[lb]{\smash{{\SetFigFont{10}{12.0}{\familydefault}{\mddefault}{\updefault}{\color[rgb]{0,0,0}...}%
}}}}
\put(9601,-4411){\makebox(0,0)[lb]{\smash{{\SetFigFont{10}{12.0}{\familydefault}{\mddefault}{\updefault}{\color[rgb]{0,0,0}...}%
}}}}
\end{picture}%

%% file: fig/convergence.tex
\begingroup
  \makeatletter
  \providecommand\color[2][]{%
    \GenericError{(gnuplot) \space\space\space\@spaces}{%
      Package color not loaded in conjunction with
      terminal option `colourtext'%
    }{See the gnuplot documentation for explanation.%
    }{Either use 'blacktext' in gnuplot or load the package
      color.sty in LaTeX.}%
    \renewcommand\color[2][]{}%
  }%
  \providecommand\includegraphics[2][]{%
    \GenericError{(gnuplot) \space\space\space\@spaces}{%
      Package graphicx or graphics not loaded%
    }{See the gnuplot documentation for explanation.%
    }{The gnuplot epslatex terminal needs graphicx.sty or graphics.sty.}%
    \renewcommand\includegraphics[2][]{}%
  }%
  \providecommand\rotatebox[2]{#2}%
  \@ifundefined{ifGPcolor}{%
    \newif\ifGPcolor
    \GPcolorfalse
  }{}%
  \@ifundefined{ifGPblacktext}{%
    \newif\ifGPblacktext
    \GPblacktexttrue
  }{}%
  \let\gplgaddtomacro\g@addto@macro
  \gdef\gplbacktext{}%
  \gdef\gplfronttext{}%
  \makeatother
  \ifGPblacktext
    \def\colorrgb#1{}%
    \def\colorgray#1{}%
  \else
    \ifGPcolor
      \def\colorrgb#1{\color[rgb]{#1}}%
      \def\colorgray#1{\color[gray]{#1}}%
      \expandafter\def\csname LTw\endcsname{\color{white}}%
      \expandafter\def\csname LTb\endcsname{\color{black}}%
      \expandafter\def\csname LTa\endcsname{\color{black}}%
      \expandafter\def\csname LT0\endcsname{\color[rgb]{1,0,0}}%
      \expandafter\def\csname LT1\endcsname{\color[rgb]{0,1,0}}%
      \expandafter\def\csname LT2\endcsname{\color[rgb]{0,0,1}}%
      \expandafter\def\csname LT3\endcsname{\color[rgb]{1,0,1}}%
      \expandafter\def\csname LT4\endcsname{\color[rgb]{0,1,1}}%
      \expandafter\def\csname LT5\endcsname{\color[rgb]{1,1,0}}%
      \expandafter\def\csname LT6\endcsname{\color[rgb]{0,0,0}}%
      \expandafter\def\csname LT7\endcsname{\color[rgb]{1,0.3,0}}%
      \expandafter\def\csname LT8\endcsname{\color[rgb]{0.5,0.5,0.5}}%
    \else
      \def\colorrgb#1{\color{black}}%
      \def\colorgray#1{\color[gray]{#1}}%
      \expandafter\def\csname LTw\endcsname{\color{white}}%
      \expandafter\def\csname LTb\endcsname{\color{black}}%
      \expandafter\def\csname LTa\endcsname{\color{black}}%
      \expandafter\def\csname LT0\endcsname{\color{black}}%
      \expandafter\def\csname LT1\endcsname{\color{black}}%
      \expandafter\def\csname LT2\endcsname{\color{black}}%
      \expandafter\def\csname LT3\endcsname{\color{black}}%
      \expandafter\def\csname LT4\endcsname{\color{black}}%
      \expandafter\def\csname LT5\endcsname{\color{black}}%
      \expandafter\def\csname LT6\endcsname{\color{black}}%
      \expandafter\def\csname LT7\endcsname{\color{black}}%
      \expandafter\def\csname LT8\endcsname{\color{black}}%
    \fi
  \fi
    \setlength{\unitlength}{0.0500bp}%
    \ifx\gptboxheight\undefined%
      \newlength{\gptboxheight}%
      \newlength{\gptboxwidth}%
      \newsavebox{\gptboxtext}%
    \fi%
    \setlength{\fboxrule}{0.5pt}%
    \setlength{\fboxsep}{1pt}%
\begin{picture}(4400.00,5280.00)%
    \gplgaddtomacro\gplbacktext{%
      \colorrgb{0.00,0.00,0.00}%
      \put(512,580){\makebox(0,0)[r]{\strut{}-6}}%
      \colorrgb{0.00,0.00,0.00}%
      \put(512,880){\makebox(0,0)[r]{\strut{}-5}}%
      \colorrgb{0.00,0.00,0.00}%
      \put(512,1180){\makebox(0,0)[r]{\strut{}-4}}%
      \colorrgb{0.00,0.00,0.00}%
      \put(512,1481){\makebox(0,0)[r]{\strut{}-3}}%
      \colorrgb{0.00,0.00,0.00}%
      \put(512,1781){\makebox(0,0)[r]{\strut{}-2}}%
      \colorrgb{0.00,0.00,0.00}%
      \put(512,2081){\makebox(0,0)[r]{\strut{}-1}}%
      \colorrgb{0.00,0.00,0.00}%
      \put(512,2381){\makebox(0,0)[r]{\strut{}0}}%
      \colorrgb{0.00,0.00,0.00}%
      \put(572,480){\makebox(0,0){\strut{}0}}%
      \colorrgb{0.00,0.00,0.00}%
      \put(1254,480){\makebox(0,0){\strut{}5}}%
      \colorrgb{0.00,0.00,0.00}%
      \put(1936,480){\makebox(0,0){\strut{}10}}%
      \colorrgb{0.00,0.00,0.00}%
      \put(2617,480){\makebox(0,0){\strut{}15}}%
      \colorrgb{0.00,0.00,0.00}%
      \put(3299,480){\makebox(0,0){\strut{}20}}%
      \colorrgb{0.00,0.00,0.00}%
      \put(3981,480){\makebox(0,0){\strut{}25}}%
    }%
    \gplgaddtomacro\gplfronttext{%
      \colorrgb{0.00,0.00,0.00}%
      \put(282,1480){\rotatebox{90}{\makebox(0,0){\strut{}$\log_{10}(\frac{1}{N}\|\mathbf{g}_k\|)$}}}%
      \colorrgb{0.00,0.00,0.00}%
      \put(2276,330){\makebox(0,0){\strut{}Time (min)}}%
      \csname LTb\endcsname%
      \put(2276,2481){\makebox(0,0){\strut{}(b)}}%
      \colorrgb{0.00,0.00,0.00}%
      \put(3678,2294){\makebox(0,0){\tiny \bf ALS\hspace{5ex}}}%
      \colorrgb{0.00,0.00,0.00}%
      \put(3678,2201){\makebox(0,0){\tiny \bf ALS-NCG\hspace{10ex}}}%
    }%
    \gplgaddtomacro\gplbacktext{%
      \colorrgb{0.00,0.00,0.00}%
      \put(512,3082){\makebox(0,0)[r]{\strut{}-6}}%
      \colorrgb{0.00,0.00,0.00}%
      \put(512,3382){\makebox(0,0)[r]{\strut{}-5}}%
      \colorrgb{0.00,0.00,0.00}%
      \put(512,3682){\makebox(0,0)[r]{\strut{}-4}}%
      \colorrgb{0.00,0.00,0.00}%
      \put(512,3983){\makebox(0,0)[r]{\strut{}-3}}%
      \colorrgb{0.00,0.00,0.00}%
      \put(512,4283){\makebox(0,0)[r]{\strut{}-2}}%
      \colorrgb{0.00,0.00,0.00}%
      \put(512,4583){\makebox(0,0)[r]{\strut{}-1}}%
      \colorrgb{0.00,0.00,0.00}%
      \put(512,4883){\makebox(0,0)[r]{\strut{}0}}%
      \colorrgb{0.00,0.00,0.00}%
      \put(572,2982){\makebox(0,0){\strut{}0}}%
      \colorrgb{0.00,0.00,0.00}%
      \put(1253,2982){\makebox(0,0){\strut{}20}}%
      \colorrgb{0.00,0.00,0.00}%
      \put(1935,2982){\makebox(0,0){\strut{}40}}%
      \colorrgb{0.00,0.00,0.00}%
      \put(2616,2982){\makebox(0,0){\strut{}60}}%
      \colorrgb{0.00,0.00,0.00}%
      \put(3298,2982){\makebox(0,0){\strut{}80}}%
      \colorrgb{0.00,0.00,0.00}%
      \put(3979,2982){\makebox(0,0){\strut{}100}}%
    }%
    \gplgaddtomacro\gplfronttext{%
      \colorrgb{0.00,0.00,0.00}%
      \put(282,3982){\rotatebox{90}{\makebox(0,0){\strut{}$\log_{10}(\frac{1}{N}\|\mathbf{g}_k\|)$}}}%
      \colorrgb{0.00,0.00,0.00}%
      \put(2276,2832){\makebox(0,0){\strut{}Iteration, $k$}}%
      \csname LTb\endcsname%
      \put(2276,4983){\makebox(0,0){\strut{}(a)}}%
      \colorrgb{0.00,0.00,0.00}%
      \put(3678,4791){\makebox(0,0){\tiny \bf ALS\hspace{5ex}}}%
      \colorrgb{0.00,0.00,0.00}%
      \put(3678,4703){\makebox(0,0){\tiny \bf  ALS-NCG\hspace{10ex}}}%
    }%
    \gplbacktext
    \put(0,0){\includegraphics{convergence}}%
    \gplfronttext
  \end{picture}%
\endgroup

%% file: fig/gran_scaling.tex
\begingroup
  \makeatletter
  \providecommand\color[2][]{%
    \GenericError{(gnuplot) \space\space\space\@spaces}{%
      Package color not loaded in conjunction with
      terminal option `colourtext'%
    }{See the gnuplot documentation for explanation.%
    }{Either use 'blacktext' in gnuplot or load the package
      color.sty in LaTeX.}%
    \renewcommand\color[2][]{}%
  }%
  \providecommand\includegraphics[2][]{%
    \GenericError{(gnuplot) \space\space\space\@spaces}{%
      Package graphicx or graphics not loaded%
    }{See the gnuplot documentation for explanation.%
    }{The gnuplot epslatex terminal needs graphicx.sty or graphics.sty.}%
    \renewcommand\includegraphics[2][]{}%
  }%
  \providecommand\rotatebox[2]{#2}%
  \@ifundefined{ifGPcolor}{%
    \newif\ifGPcolor
    \GPcolorfalse
  }{}%
  \@ifundefined{ifGPblacktext}{%
    \newif\ifGPblacktext
    \GPblacktexttrue
  }{}%
  \let\gplgaddtomacro\g@addto@macro
  \gdef\gplbacktext{}%
  \gdef\gplfronttext{}%
  \makeatother
  \ifGPblacktext
    \def\colorrgb#1{}%
    \def\colorgray#1{}%
  \else
    \ifGPcolor
      \def\colorrgb#1{\color[rgb]{#1}}%
      \def\colorgray#1{\color[gray]{#1}}%
      \expandafter\def\csname LTw\endcsname{\color{white}}%
      \expandafter\def\csname LTb\endcsname{\color{black}}%
      \expandafter\def\csname LTa\endcsname{\color{black}}%
      \expandafter\def\csname LT0\endcsname{\color[rgb]{1,0,0}}%
      \expandafter\def\csname LT1\endcsname{\color[rgb]{0,1,0}}%
      \expandafter\def\csname LT2\endcsname{\color[rgb]{0,0,1}}%
      \expandafter\def\csname LT3\endcsname{\color[rgb]{1,0,1}}%
      \expandafter\def\csname LT4\endcsname{\color[rgb]{0,1,1}}%
      \expandafter\def\csname LT5\endcsname{\color[rgb]{1,1,0}}%
      \expandafter\def\csname LT6\endcsname{\color[rgb]{0,0,0}}%
      \expandafter\def\csname LT7\endcsname{\color[rgb]{1,0.3,0}}%
      \expandafter\def\csname LT8\endcsname{\color[rgb]{0.5,0.5,0.5}}%
    \else
      \def\colorrgb#1{\color{black}}%
      \def\colorgray#1{\color[gray]{#1}}%
      \expandafter\def\csname LTw\endcsname{\color{white}}%
      \expandafter\def\csname LTb\endcsname{\color{black}}%
      \expandafter\def\csname LTa\endcsname{\color{black}}%
      \expandafter\def\csname LT0\endcsname{\color{black}}%
      \expandafter\def\csname LT1\endcsname{\color{black}}%
      \expandafter\def\csname LT2\endcsname{\color{black}}%
      \expandafter\def\csname LT3\endcsname{\color{black}}%
      \expandafter\def\csname LT4\endcsname{\color{black}}%
      \expandafter\def\csname LT5\endcsname{\color{black}}%
      \expandafter\def\csname LT6\endcsname{\color{black}}%
      \expandafter\def\csname LT7\endcsname{\color{black}}%
      \expandafter\def\csname LT8\endcsname{\color{black}}%
    \fi
  \fi
    \setlength{\unitlength}{0.0500bp}%
    \ifx\gptboxheight\undefined%
      \newlength{\gptboxheight}%
      \newlength{\gptboxwidth}%
      \newsavebox{\gptboxtext}%
    \fi%
    \setlength{\fboxrule}{0.5pt}%
    \setlength{\fboxsep}{1pt}%
\begin{picture}(5200.00,3900.00)%
    \gplgaddtomacro\gplbacktext{%
      \colorrgb{0.00,0.00,0.00}%
      \put(888,768){\makebox(0,0)[r]{\strut{}0}}%
      \colorrgb{0.00,0.00,0.00}%
      \put(888,1479){\makebox(0,0)[r]{\strut{}50}}%
      \colorrgb{0.00,0.00,0.00}%
      \put(888,2190){\makebox(0,0)[r]{\strut{}100}}%
      \colorrgb{0.00,0.00,0.00}%
      \put(888,2900){\makebox(0,0)[r]{\strut{}150}}%
      \colorrgb{0.00,0.00,0.00}%
      \put(888,3611){\makebox(0,0)[r]{\strut{}200}}%
      \colorrgb{0.00,0.00,0.00}%
      \put(1343,528){\makebox(0,0){\strut{}1}}%
      \colorrgb{0.00,0.00,0.00}%
      \put(1966,528){\makebox(0,0){\strut{}2}}%
      \colorrgb{0.00,0.00,0.00}%
      \put(2588,528){\makebox(0,0){\strut{}3}}%
      \colorrgb{0.00,0.00,0.00}%
      \put(3211,528){\makebox(0,0){\strut{}4}}%
      \colorrgb{0.00,0.00,0.00}%
      \put(3833,528){\makebox(0,0){\strut{}5}}%
      \colorrgb{0.00,0.00,0.00}%
      \put(4456,528){\makebox(0,0){\strut{}6}}%
    }%
    \gplgaddtomacro\gplfronttext{%
      \colorrgb{0.00,0.00,0.00}%
      \put(192,2189){\rotatebox{90}{\makebox(0,0){\strut{}Time Per Iteration (s)}}}%
      \colorrgb{0.00,0.00,0.00}%
      \put(2899,168){\makebox(0,0){\strut{}$n_u$ (millions)}}%
      \colorrgb{0.00,0.00,0.00}%
      \put(4128,3488){\makebox(0,0)[r]{\strut{}\footnotesize\bf ALS}}%
      \colorrgb{0.00,0.00,0.00}%
      \put(4128,3368){\makebox(0,0)[r]{\strut{}\footnotesize\bf ALS-NCG}}%
    }%
    \gplbacktext
    \put(0,0){\includegraphics{gran_scaling}}%
    \gplfronttext
  \end{picture}%
\endgroup

%% file: fig/speedup.tex
\begingroup
  \makeatletter
  \providecommand\color[2][]{%
    \GenericError{(gnuplot) \space\space\space\@spaces}{%
      Package color not loaded in conjunction with
      terminal option `colourtext'%
    }{See the gnuplot documentation for explanation.%
    }{Either use 'blacktext' in gnuplot or load the package
      color.sty in LaTeX.}%
    \renewcommand\color[2][]{}%
  }%
  \providecommand\includegraphics[2][]{%
    \GenericError{(gnuplot) \space\space\space\@spaces}{%
      Package graphicx or graphics not loaded%
    }{See the gnuplot documentation for explanation.%
    }{The gnuplot epslatex terminal needs graphicx.sty or graphics.sty.}%
    \renewcommand\includegraphics[2][]{}%
  }%
  \providecommand\rotatebox[2]{#2}%
  \@ifundefined{ifGPcolor}{%
    \newif\ifGPcolor
    \GPcolorfalse
  }{}%
  \@ifundefined{ifGPblacktext}{%
    \newif\ifGPblacktext
    \GPblacktexttrue
  }{}%
  \let\gplgaddtomacro\g@addto@macro
  \gdef\gplbacktext{}%
  \gdef\gplfronttext{}%
  \makeatother
  \ifGPblacktext
    \def\colorrgb#1{}%
    \def\colorgray#1{}%
  \else
    \ifGPcolor
      \def\colorrgb#1{\color[rgb]{#1}}%
      \def\colorgray#1{\color[gray]{#1}}%
      \expandafter\def\csname LTw\endcsname{\color{white}}%
      \expandafter\def\csname LTb\endcsname{\color{black}}%
      \expandafter\def\csname LTa\endcsname{\color{black}}%
      \expandafter\def\csname LT0\endcsname{\color[rgb]{1,0,0}}%
      \expandafter\def\csname LT1\endcsname{\color[rgb]{0,1,0}}%
      \expandafter\def\csname LT2\endcsname{\color[rgb]{0,0,1}}%
      \expandafter\def\csname LT3\endcsname{\color[rgb]{1,0,1}}%
      \expandafter\def\csname LT4\endcsname{\color[rgb]{0,1,1}}%
      \expandafter\def\csname LT5\endcsname{\color[rgb]{1,1,0}}%
      \expandafter\def\csname LT6\endcsname{\color[rgb]{0,0,0}}%
      \expandafter\def\csname LT7\endcsname{\color[rgb]{1,0.3,0}}%
      \expandafter\def\csname LT8\endcsname{\color[rgb]{0.5,0.5,0.5}}%
    \else
      \def\colorrgb#1{\color{black}}%
      \def\colorgray#1{\color[gray]{#1}}%
      \expandafter\def\csname LTw\endcsname{\color{white}}%
      \expandafter\def\csname LTb\endcsname{\color{black}}%
      \expandafter\def\csname LTa\endcsname{\color{black}}%
      \expandafter\def\csname LT0\endcsname{\color{black}}%
      \expandafter\def\csname LT1\endcsname{\color{black}}%
      \expandafter\def\csname LT2\endcsname{\color{black}}%
      \expandafter\def\csname LT3\endcsname{\color{black}}%
      \expandafter\def\csname LT4\endcsname{\color{black}}%
      \expandafter\def\csname LT5\endcsname{\color{black}}%
      \expandafter\def\csname LT6\endcsname{\color{black}}%
      \expandafter\def\csname LT7\endcsname{\color{black}}%
      \expandafter\def\csname LT8\endcsname{\color{black}}%
    \fi
  \fi
    \setlength{\unitlength}{0.0500bp}%
    \ifx\gptboxheight\undefined%
      \newlength{\gptboxheight}%
      \newlength{\gptboxwidth}%
      \newsavebox{\gptboxtext}%
    \fi%
    \setlength{\fboxrule}{0.5pt}%
    \setlength{\fboxsep}{1pt}%
\begin{picture}(5200.00,3900.00)%
    \gplgaddtomacro\gplbacktext{%
      \colorrgb{0.00,0.00,0.00}%
      \put(650,832){\makebox(0,0)[r]{\strut{}0}}%
      \colorrgb{0.00,0.00,0.00}%
      \put(650,1291){\makebox(0,0)[r]{\strut{}1}}%
      \colorrgb{0.00,0.00,0.00}%
      \put(650,1750){\makebox(0,0)[r]{\strut{}2}}%
      \colorrgb{0.00,0.00,0.00}%
      \put(650,2210){\makebox(0,0)[r]{\strut{}3}}%
      \colorrgb{0.00,0.00,0.00}%
      \put(650,2669){\makebox(0,0)[r]{\strut{}4}}%
      \colorrgb{0.00,0.00,0.00}%
      \put(650,3128){\makebox(0,0)[r]{\strut{}5}}%
      \colorrgb{0.00,0.00,0.00}%
      \put(650,3587){\makebox(0,0)[r]{\strut{}6}}%
      \colorrgb{0.00,0.00,0.00}%
      \put(1119,572){\makebox(0,0){\strut{}$10^{-3}$}}%
      \colorrgb{0.00,0.00,0.00}%
      \put(2158,572){\makebox(0,0){\strut{}$10^{-2}$}}%
      \colorrgb{0.00,0.00,0.00}%
      \put(3196,572){\makebox(0,0){\strut{}$10^{-1}$}}%
      \colorrgb{0.00,0.00,0.00}%
      \put(4235,572){\makebox(0,0){\strut{}$10^{0}$}}%
    }%
    \gplgaddtomacro\gplfronttext{%
      \colorrgb{0.00,0.00,0.00}%
      \put(208,2209){\rotatebox{90}{\makebox(0,0){\strut{}Relative Speedup}}}%
      \colorrgb{0.00,0.00,0.00}%
      \put(2768,182){\makebox(0,0){\strut{}$\frac{1}{N}\|\mathbf{g}_k\|$}}%
      \colorrgb{0.00,0.00,0.00}%
      \put(4044,3459){\makebox(0,0)[r]{\strut{}\bf\footnotesize 1M users}}%
      \colorrgb{0.00,0.00,0.00}%
      \put(4044,3329){\makebox(0,0)[r]{\strut{}\bf\footnotesize 3M users}}%
      \colorrgb{0.00,0.00,0.00}%
      \put(4044,3199){\makebox(0,0)[r]{\strut{}\bf\footnotesize 6M users}}%
    }%
    \gplbacktext
    \put(0,0){\includegraphics{speedup}}%
    \gplfronttext
  \end{picture}%
\endgroup

%% file: main.bbl
\begin{thebibliography}{10}
\providecommand{\url}[1]{#1}
\csname url@samestyle\endcsname
\providecommand{\newblock}{\relax}
\providecommand{\bibinfo}[2]{#2}
\providecommand{\BIBentrySTDinterwordspacing}{\spaceskip=0pt\relax}
\providecommand{\BIBentryALTinterwordstretchfactor}{4}
\providecommand{\BIBentryALTinterwordspacing}{\spaceskip=\fontdimen2\font plus
\BIBentryALTinterwordstretchfactor\fontdimen3\font minus
  \fontdimen4\font\relax}
\providecommand{\BIBforeignlanguage}[2]{{%
\expandafter\ifx\csname l@#1\endcsname\relax
\typeout{** WARNING: IEEEtran.bst: No hyphenation pattern has been}%
\typeout{** loaded for the language `#1'. Using the pattern for}%
\typeout{** the default language instead.}%
\else
\language=\csname l@#1\endcsname
\fi
#2}}
\providecommand{\BIBdecl}{\relax}
\BIBdecl

\bibitem{Sarwar:2001}
B.~Sarwar, G.~Karypis, J.~Konstan, and J.~Riedl, ``Item-based collaborative
  filtering recommendation algorithms,'' in \emph{Proceedings of the 10th
  International Conference on World Wide Web}, 2001, pp. 285--295.

\bibitem{koren2008factorization}
Y.~Koren, ``Factorization meets the neighborhood: a multifaceted collaborative
  filtering model,'' in \emph{14th ACM SIGKDD international conference on
  Knowledge discovery and data mining}.\hskip 1em plus 0.5em minus 0.4em\relax
  ACM, 2008, pp. 426--434.

\bibitem{Bobadilla:2013}
J.~Bobadilla, F.~Ortega, A.~Hernando, and A.~Guti{\'e}Rrez, ``Recommender
  systems survey,'' \emph{Knowledge-Based Systems}, vol.~46, pp. 109--132,
  2013.

\bibitem{Linden:2003}
G.~Linden, B.~Smith, and J.~York, ``Amazon.com recommendations: Item-to-item
  collaborative filtering,'' \emph{IEEE Internet Comput.}, vol.~7, no.~1, pp.
  76--80, 2003.

\bibitem{Bell:2007a}
R.~M. Bell and Y.~Koren, ``Lessons from the {N}etflix prize challenge,''
  \emph{SIGKDD Explor. Newsl.}, vol.~9, no.~2, pp. 75--79, 2007.

\bibitem{Johnson:2014}
C.~C. Johnson, ``Logistic matrix factorization for implicit feedback data,'' in
  \emph{NIPS Workshop on Distributed Machine Learning and Matrix Computations},
  2014.

\bibitem{Dror:2012}
G.~Dror, N.~Koenigstein, Y.~Koren, and M.~Weimer, ``The {Y}ahoo! music dataset
  and {KDD}-cup '11,'' \emph{JMLR: Workshop and Conference Proceedings},
  vol.~18, pp. 3--18, 2012.

\bibitem{Koren:2009}
Y.~Koren, R.~Bell, and C.~Volinsky, ``Matrix factorization techniques for
  recommender systems,'' \emph{Computer}, vol.~42, no.~8, pp. 30--37, 2009.

\bibitem{Funk:2006}
S.~Funk, ``Netflix update: Try this at home,''
  http://sifter.org/~simon/journal/20061211.html, 2006.

\bibitem{Hu:2008}
Y.~Hu, Y.~Koren, and C.~Volinsky, ``Collaborative filtering for implicit
  feedback datasets,'' in \emph{Proceedings of the 2008 Eighth IEEE
  International Conference on Data Mining}, 2008, pp. 263--272.

\bibitem{Nocedal:2006}
J.~Nocedal and S.~J. Wright, \emph{Numerical Optimization}, 2nd~ed.\hskip 1em
  plus 0.5em minus 0.4em\relax New York: Springer, 2006.

\bibitem{DeSterck:2015}
H.~{De Sterck} and M.~Winlaw, ``A nonlinearly preconditioned conjugate gradient
  algorithm for rank-{$R$} canonical tensor approximation,'' \emph{Numer.
  Linear Algebra Appl.}, vol.~22, pp. 410--432, 2015.

\bibitem{PETSC:2015}
P.~Brune, M.~G. Knepley, B.~F. Smith, and X.~Tu, ``Composing scalable nonlinear
  algebraic solvers,'' \emph{SIAM Review}, forthcoming.

\bibitem{DeSterck:2012}
H.~{De Sterck}, ``A nonlinear {GMRES} optimization algorithm for canonical
  tensor decomposition,'' \emph{{SIAM} J. Sci. Comput.}, vol.~34, pp.
  A1351--A1379, 2012.

\bibitem{Fang:2009}
H.~Fang and Y.~Saad, ``Two classes of multisecant methods for nonlinear
  acceleration,'' \emph{Numer. Linear Algebra Appl.}, vol.~16, pp. 197--221,
  2009.

\bibitem{Walker:2011}
H.~Walker and P.~Ni, ``Anderson acceleration for fixed-point iterations,''
  \emph{{SIAM} J. Numer. Anal.}, vol.~49, pp. 1715--1735, 2011.

\bibitem{Anderson:1965}
D.~Anderson, ``Iterative procedures for nonlinear integral equations,''
  \emph{J. ACM}, vol.~12, pp. 547--560, 1965.

\bibitem{Concus:1977}
P.~Concus, G.~H. Golub, and D.~P. O'Leary, ``Numerical solution of nonlinear
  elliptical partial differential equations by a generalized conjugate gradient
  method,'' \emph{Computing}, vol.~19, pp. 321--339, 1977.

\bibitem{gemulla2011large}
R.~Gemulla, E.~Nijkamp, P.~J. Haas, and Y.~Sismanis, ``Large-scale matrix
  factorization with distributed stochastic gradient descent,'' in
  \emph{Proceedings of the 17th ACM SIGKDD international conference on
  Knowledge discovery and data mining}.\hskip 1em plus 0.5em minus 0.4em\relax
  ACM, 2011, pp. 69--77.

\bibitem{teflioudi2012distributed}
C.~Teflioudi, F.~Makari, and R.~Gemulla, ``Distributed matrix completion,'' in
  \emph{12th International Conference on Data Mining (ICDM)}.\hskip 1em plus
  0.5em minus 0.4em\relax IEEE, 2012, pp. 655--664.

\bibitem{yu2012scalable}
H.-F. Yu, C.-J. Hsieh, I.~Dhillon \emph{et~al.}, ``Scalable coordinate descent
  approaches to parallel matrix factorization for recommender systems,'' in
  \emph{12th International Conference on Data Mining (ICDM)}.\hskip 1em plus
  0.5em minus 0.4em\relax IEEE, 2012, pp. 765--774.

\bibitem{zaharia2012resilient}
M.~Zaharia, M.~Chowdhury, T.~Das, A.~Dave, J.~Ma, M.~McCauley, M.~J. Franklin,
  S.~Shenker, and I.~Stoica, ``Resilient distributed datasets: {A}
  fault-tolerant abstraction for in-memory cluster computing,'' in
  \emph{Proceedings of the 9th USENIX conference on Networked Systems Design
  and Implementation}, 2012, pp. 15--28.

\bibitem{Zhou:2008}
Y.~Zhou, D.~Wilkinson, R.~Schreiber, and R.~Pan, ``Large-scale parallel
  collaborative filtering for the {N}etflix prize,'' in \emph{Proceedings of
  the 4th International Conference on Algorithmic Aspects in Information and
  Management}, 2008, pp. 337--348.

\bibitem{Johnson2014}
\BIBentryALTinterwordspacing
C.~Johnson, ``Music recommendations at scale with {Spark},'' 2014, {Spark
  Summit}. [Online]. Available:
  \url{https://spark-summit.org/2014/talk/music-recommendations-at-scale-with-spark}
\BIBentrySTDinterwordspacing

\bibitem{MovieLens:2015}
\BIBentryALTinterwordspacing
``{M}ovie{L}ens {20M} dataset,'' March 2015, {G}roup{L}ens. [Online].
  Available: \url{http://grouplens.org/datasets/movielens/20m/}
\BIBentrySTDinterwordspacing

\bibitem{Tikhonov:77}
A.~N. Tikhonov and V.~Y. Arsenin, \emph{Solutions of ill-posed problems}.\hskip
  1em plus 0.5em minus 0.4em\relax New York: John Wiley, 1977.

\bibitem{Polak:1969}
E.~Polak and G.~Ribi\`{e}re, ``Note sur la convergence de m\'{e}thodes de
  directions conjug\'{e}es,'' \emph{Revue Fran\c{c}asie d`Informatique et de
  Recherche Op\'{e}rationnelle}, vol.~16, pp. 35--43, 1969.

\bibitem{Poblano}
D.~M. Dunlavy, T.~G. Kolda, and E.~Acar, ``Poblano v1.0: {A} {MATLAB} toolbox
  for gradient-based optimization,'' Sandia National Laboratories, Albuquerque,
  NM and Livermore, CA, Tech. Rep. SAND2010-1422, March 2010.

\bibitem{Kendall:1938}
M.~G. Kendall, ``A new measure of rank correlation,'' \emph{Biometrika},
  vol.~30, no. 1/2, pp. 81--93, 1938.

\bibitem{zaharia2010delay}
M.~Zaharia, D.~Borthakur, J.~Sen~Sarma, K.~Elmeleegy, S.~Shenker, and
  I.~Stoica, ``Delay scheduling: {A} simple technique for achieving locality
  and fairness in cluster scheduling,'' in \emph{Proceedings of the 5th
  European Conference on Computer systems}, 2010, pp. 265--278.

\bibitem{blackford2002updated}
L.~S. Blackford, A.~Petitet, R.~Pozo, K.~Remington, R.~C. Whaley, J.~Demmel,
  J.~Dongarra, I.~Duff, S.~Hammarling, G.~Henry \emph{et~al.}, ``An updated set
  of basic linear algebra subprograms {(BLAS)},'' \emph{ACM Trans. Math.
  Software}, vol.~28, no.~2, pp. 135--151, 2002.

\bibitem{anderson1999lapack}
E.~Anderson, Z.~Bai, C.~Bischof, S.~Blackford, J.~Demmel, J.~Dongarra,
  J.~Du~Croz, A.~Greenbaum, S.~Hammerling, A.~McKenney \emph{et~al.},
  \emph{LAPACK Users' Guide}.\hskip 1em plus 0.5em minus 0.4em\relax {SIAM},
  1999, vol.~9.

\bibitem{davidson2013optimizing}
A.~Davidson and A.~Or, ``Optimizing shuffle performance in {Spark},''
  \emph{University of California, Berkeley-Department of Electrical Engineering
  and Computer Sciences, Tech. Rep}, 2013.

\bibitem{Dai2015}
\BIBentryALTinterwordspacing
J.~Dai, ``Experience and lessons learned for large-scale graph analysis using
  {GraphX},'' 2015, {Spark Summit East}. [Online]. Available:
  \url{https://spark-summit.org/east-2015/talk/experience-and-lessons-learned-for-large-scale-graph-analysis-using-graphx}
\BIBentrySTDinterwordspacing

\end{thebibliography}
